%% file: main.tex
\documentclass[letterpaper, 10 pt, journal]{IEEEtran}
\usepackage{cite}
\usepackage{amsmath} 
\usepackage{amssymb}  
\usepackage{amsfonts}
\usepackage{mathrsfs}
\usepackage{mathtools}
\usepackage[dvipsnames]{xcolor}
\usepackage[thmmarks, amsmath]{ntheorem}
\usepackage{enumerate}
\usepackage{circuitikz}
\usepackage{pgfplots}
\usepackage{booktabs}
\usepackage[bbgreekl]{mathbbol}
\usepackage{colortbl}
\usepackage{ragged2e}
\usepackage{enumitem}
\usepackage{relsize}
\usepackage{standalone}
\usepackage{tikz}
\usepackage[ruled,vlined]{algorithm2e}
\usepackage{etoolbox}

\makeatletter
\patchcmd{\@algocf@start}
  {-1.5em}
  {0pt}
  {}{}
\makeatother
\renewcommand{\baselinestretch}{1.0}
\usetikzlibrary{positioning, shapes}
\usetikzlibrary{decorations.pathreplacing}
\usetikzlibrary{arrows}
\usetikzlibrary{plotmarks}
\usetikzlibrary{patterns}
\usetikzlibrary{spy}
\usetikzlibrary{calc,trees,positioning,arrows,chains,shapes.geometric,%
    decorations.pathreplacing,decorations.pathmorphing,shapes,decorations.markings,matrix,shapes.symbols}

\tikzset{
>=stealth',
  punktchain/.style={rectangle, rounded corners, 
    draw=black, very thick,text width=10em, 
    minimum height=3em, text centered, on chain},
  line/.style={draw, thick, <-},
  element/.style={tape,top color=white,bottom color=blue!50!black!60!,
    minimum width=8em,draw=blue!40!black!90, very thick,
    text width=10em, minimum height=3.5em, text centered, on chain},
  every join/.style={->, thick,shorten >=1pt},
  tuborg/.style={decorate},
  tubnode/.style={midway, right=2pt},
}

\graphicspath{{./figures/}}

\usepackage[font=small,skip=2pt]{caption}

\newcommand{\ignore}[1]{}
\allowdisplaybreaks
\DeclareSymbolFont{bbold}{U}{bbold}{m}{n}
\DeclareSymbolFontAlphabet{\mathbbold}{bbold}

\newcommand{\oprocendsymbol}{\hbox{$\bullet$}}
\newcommand{\oprocend}{\relax\ifmmode\else\unskip\hfill\fi\oprocendsymbol}
\newtheorem{theorem}{Theorem}
\newtheorem{lemma}[theorem]{Lemma}

\newtheorem{remark}[theorem]{Remark}

\usepackage{tikz}
\usetikzlibrary{external}
\usepackage{standalone}

\newcolumntype{R}[1]{>{\RaggedLeft\arraybackslash}p{#1}}

\usepackage[absolute,overlay]{textpos}
\usepackage{pgfplots}
\pgfplotsset{compat=1.6,
        scaled x ticks = false, 
        xticklabel style={/pgf/number format/fixed,/pgf/number format/precision=3},
        legend image post style={xscale=0.3},
             legend image post style={yscale=0.6}
        } 

\newcolumntype{P}[1]{>{\centering\arraybackslash}p{#1}}
\definecolor{mycoloryell}{rgb}{0.92941,0.69412,0.12549}%
\definecolor{mycolorblue}{rgb}{0.30588,0.39608,0.58039}%
\usepackage[ruled,vlined]{algorithm2e}

\begin{document}

\begin{textblock*}{\textwidth}(15mm,9mm) 
\centering \bf \textcolor{NavyBlue}{To appear in the {\emph{IEEE Transactions on Smart Grid}}\\ https://doi.org/10.1109/TSG.2021.3122879} \\
\end{textblock*}

\title{Pricing and Energy Trading in Peer-to-peer \\ Zero Marginal-cost Microgrids}
\author{Jonathan Lee, Rodrigo Henriquez-Auba, Bala Kameshwar Poolla, and Duncan S. Callaway
\thanks{Jonathan Lee and Duncan S. Callaway are with the Energy and Resources Group, and Rodrigo Henriquez-Auba is with the Department of Electrical Engineering and Computer Sciences, University of California Berkeley, CA. Email: {\tt \{jtlee, dcal, rhenriquez\}@berkeley.edu}}%
\thanks{Bala Kameshwar Poolla is with the National Renewable Energy Laboratory, Golden, CO, USA. Email: {\tt bpoolla@nrel.gov}}%
\thanks{This research is supported by the National Science Foundation under grants CPS-1646612 and CyberSEES-1539585 and the Graduate Research Fellowship Program.}
\thanks{This work was authored in part by the National Renewable Energy Laboratory (NREL), operated by Alliance for Sustainable Energy, LLC, for the U.S. Department of Energy (DOE) under Contract DE-AC36-08GO28308 and supported by the Laboratory Directed Research and Development (LDRD) Program at NREL. The views expressed in the paper do not necessarily represent the views of the DOE or the U.S. Government. The U.S. Government retains and the publisher, by accepting the paper for publication, acknowledges that the U.S. Government retains a nonexclusive, paid-up, irrevocable, worldwide license to publish or reproduce the published form of this work, or allow others to do so, for U.S. Government purposes.}
}
\maketitle

\begin{abstract}
Efforts to utilize $100$\% renewable energy in community microgrids require new approaches to energy markets and transactions to efficiently address periods of scarce energy supply. In this paper we contribute to the promising approach of peer-to-peer (P2P) energy trading in two main ways: analysis of a centralized, welfare-maximizing economic dispatch that characterizes optimal price and allocations, and a novel P2P system for negotiating energy trades that yields physically feasible and at least weakly Pareto-optimal outcomes. Our main results are 1) that optimal pricing is insufficient to induce agents with batteries to take optimal actions, 2) a novel P2P algorithm to address this while keeping private information, 3) a formal proof that this algorithm converges to the centralized solution in the case of two agents negotiating for a single period, and  4) numerical simulations of the P2P algorithm performance with up to $10$ agents and $24$ periods that show it converges on average to total welfare within $0.1$\% of the social optimum in on the order of $10$s to $100$s of iterations, increasing with the number of agents, time periods, and total storage capacity.
\end{abstract}


\section*{Nomenclature}
\noindent\hspace{-0.25cm}\textit{Sets and indices}
\begin{IEEEdescription}[\IEEEusemathlabelsep\IEEEsetlabelwidth{$P_{11}$}]
\item[$\mathcal{C}$] Set of agents/consumers indexed by $n$.
\item[$\mathcal{U}$] Subset of agents that propose quantities ($q$-agents), indexed by $k$.
\item[$\mathcal{X}$] Subset of $q$-agents that have exited the negotiation.
\item[$\mathcal{Y}$] Subset of $q$-agents that are still negotiating.
\item[$\mathcal{V}$] Subset of agents that respond with price ($\pi$-agents), indexed by $v$.
\item[$\mathcal{T}$] Set of time periods, indexed by $t$.
\item[$\mathcal{B}$] Set of batteries, indexed by $i$.
\item[$\mathcal{G}$] Set of generators, indexed by $g$.
\end{IEEEdescription}
\vspace{0.4cm}
\noindent\hspace{-0.25cm}\textit{Variables}
\begin{IEEEdescription}[\IEEEusemathlabelsep\IEEEsetlabelwidth{$P_{11}$}]
\item[$d_{n,t}$] Local power consumption of agent $n$ at time $t$.
\item[$p_{g,t}^s$] PV power production from generation $g$ at time $t$.
\item[$p_{i,t}^b$] Discharge power from battery $i$ at time $t$.
\item[$s_{i,t}$] State of charge of battery $i$ at time $t$.
\item[$\lambda_{x,t}^y$] Dual variable of constraint $y$, where $x \in \{n,i,g\}$.
\item[$q_{k,t}$] Proposed quantity by $q$-agent $k$ to receive from the $\pi$-agent at time $t$.
\item[$q_{k,t}'$] Projection of $q_{k,t}$ to the feasible quantity set by the $\pi$-agent at time $t$.
\item[$\beta$] Auxiliary variable used to project $q_{k,t} \to q_{k,t}'$.
\item[$\alpha_n$] Binary state, true iff agent $n$ prefers current offer to no trade at all.
\item[$o_{k,t}$] Binary state, true iff a proposed quantity from $q$-agent $k$ at time $t$ is ``oscillating'', i.e., not monotonically increasing or decreasing over $3$ iterations.
\item[$\eta_k$] Binary state, true iff $q$-agent $k$ requests to settle.
\end{IEEEdescription}
\vspace{0.4cm}
\noindent\hspace{-0.25cm}\textit{Parameters}
\begin{IEEEdescription}[\IEEEusemathlabelsep\IEEEsetlabelwidth{$P_{11}$}]
\item[$\bar{P}_{g,t}^s$] Max. PV power from generation $g$ at time $t$.
\item[$\bar{P}_{i,t}^b$] Max. rate of charge/discharge of battery $i$ at time $t$.
\item[$\bar{S}_{i,t}^b$] Max. energy capacity of battery $i$ at time $t$.
\item[$\delta^{(0)}_k$] Initial size of step-limiting constraint on a $q$-agent $k$.
\item[$\gamma$] Shrinking rate of step-limiting constraint $\!\in\! (0,1)$.
\end{IEEEdescription}
\textit{Notation:} The utility functions of the agents with respect to their local demand are denoted by $U_n(d_n)$, marginal utility $\partial U_n/ \partial d_n$ by, $g_n(d_n)$, and its inverse $g_n^{-1}(d_n)$:=$h_n(\pi_n)$. Bold symbols represent a vector or a collection of points, e.g., $\boldsymbol{q}_k \equiv \{q_{k,t}\}_{t\in \mathcal{T}}$. The power and energy units are kW and kWh, and $\Delta T$ is the time step duration in hours. The symbols $\neg$, $\vee$, and $\wedge$ denote logical negation, OR, and AND respectively.

\section{Introduction}
Due to declining technology costs and a drive to reduce carbon emissions, $100$\% renewable electricity grids systems are receiving increasing attention. California's $2018$ Senate Bill $100$, for example, sets a large-scale $100$\% renewables target for $2045$. At community scales, $100$\% renewable microgrids for resilience and energy access in rural areas have become competitive with hybrid solutions with fuel-based generators, and can be preferable in cases where emissions, fuel logistics, or generator maintenance are strong concerns.

Novel pricing mechanisms for $100$\% renewable systems are not yet well-developed, but we contend they will become increasingly important for policy-makers and practitioners. For example, extending the current paradigm where load serving entities procure electricity at the lowest cost to meet inflexible demand, implies that zero (short-run) energy costs would lead to a zero (short-run) price \cite{strbac2021decarbonization}. However, Fripp \textit{et al} \cite{fripp2018variable} have shown by including demand-side participation in a capacity expansion model for Hawaii that dynamic electricity pricing to consumers results in non-zero prices and is increasingly important for maximizing welfare in $100$\% renewable systems. 

In this paper we contribute a general theoretical analysis of pricing in $100$\% renewable plus storage systems, characterizing optimal price dynamics and the challenges energy storage presents for standard bidding mechanisms. We focus on zero marginal cost renewables, with a solar generation case in simulations.We then propose an approach for community microgrids where individual ``prosumers'' with solar and storage could interact informally in a ``peer-to-peer'' (P2P) system to negotiate energy trades and form dynamic electricity prices and allocations that approximate optimal outcomes.

In the next two paragraphs we briefly review the relevant literature in this space, arguing specifically that the case of peer-to-peer trading in $100$\% renewable plus storage systems requires more analysis and innovation. P2P systems are valuable for grid resiliency, renewables integration, electricity access in less developed regions, and individual participation in electricity systems \cite{koepke2016against, green2017citizen, sousa2019peer, tushar2020peer,tushar2021peer}. Early work proposed centrally coordinated energy trading between distributed energy resources (DERs) where the generation and battery storage are fully controllable \cite{majumder2014efficient, cherukuri2019iterative}. In \cite{mengelkamp2018designing}, the authors lay the foundation for defining the physical and virtual layers required for a pooling-based system, but the paper does not develop bidding strategies for agents and assumes the microgrid remains connected to the main grid; \cite{oh2020peer} describes stochastic P2P methods to match prosumers with consumers and share profit, however no storage is considered in the model; similarly, \cite{paudel2019pricing, nguyen2020optimal, umer2021novel, ullah2021peer} propose ADMM-based methods to determine dispatch and/or pricing in P2P settings without storage; \cite{azim2021cooperative} proposes a cooperative coalition scheme, based on energy reduction that can achieve savings for their participants, but the method requires a pricing scheme that varies if the demand surpasses a given threshold to encourage power reduction and;\cite{kim2019p2p} describes a P2P architecture accounting for network charges, but does not consider the temporal aspects of storage, while \cite{zhang2019energy} uses comfort constraints for the next time step  as a limited approach to address this. 

Energy storage is fundamental to 100\% renewable systems, and some papers do incorporate it into P2P algorithms. For example, \cite{paudel2018peer} proposes a game-theoretic model, while a few papers define specific rules for battery charge/discharge cycles based on the traded quantity at each time step \cite{luth2018local, long2018peer}. An additional group of algorithms consider storage and either exchange shadow prices, employ ADMM-based bilateral trading mechanisms or additional cost-sharing methods  \cite{grzanic2021electricity, wang2019shadow, long2017peer, van2020integrated, alam2019peer, lyu2021fully}. However, none of these address the context of scarce, zero marginal-cost renewables coupled with storage.

We address the gap in the literature with a novel P2P approach that can describe informal interactions between prosumers negotiating trades for electricity in a finite time horizon setting. We assume P2P agents individually derive private utility from energy use, and iterate on price and quantity bids with their peers until convergence. In contrast to a centralized approach, the P2P approach maintains the privacy of individual utility functions and addresses the complexity of bidding storage while converging to the centralized solution in special cases.
The main contributions are:
\begin{enumerate}[label={(\arabic*)}, leftmargin=0.25cm, itemindent=0.3cm]
\setlength{\itemsep}{2pt}
    \item In Section II, to characterize the optimal price and explain the challenge of coordinating storage through centralized pricing, we formulate a centralized optimization model and highlight several non-trivial observations, such as optimal prices not uniquely determining battery dispatch decisions.
    \item In Section III, we define a novel P2P system with minimal prescriptive rules through which agents with private information exchange offers to arrive at a trade, and theoretically prove convergence for $2$-agents with a single time period. In the proposed algorithm, each agent can easily manage storage and state of charge constraints in its private decision.
    \item In Section IV, we find that the P2P approach converges to a solution in all of $1200$ general cases that were simulated, with welfare outcomes on average within $0.1$\% of the centralized while maintaining the privacy of utility functions, with a worst-case divergence of $8$\% that can arise from longer time horizons and relatively large storage capacities. 
\end{enumerate}

\section{Centralized Welfare Maximization Approach}
In this section, we define a model for optimal energy dispatch over a finite time horizon, analyze the solution for relevant insights into P2P electricity markets, and illustrate its dependence on energy storage through example. The model applies a utility maximization framework. For analysis of the optimal dispatch, we take the perspective of a benevolent central operator and assume knowledge of the individual utility functions. In practice this could be the perspective of a DER aggregator or a distribution system operator; however, it is difficult to know utility functions in practice, and this issue is a fundamental motivation for exploring peer-to-peer markets in the first place. We use a deterministic approach, where decisions are made off of a single, expected forecast of solar generation without hedging for uncertainty, and note this can be suboptimal to stochastic approaches. We also present an idealized battery model for simplicity, but show in Appendix~\ref{sec:detailed_battery_model} that the key insights still hold when we incorporate constant charge and discharge inefficiencies, self-discharge, and asymmetric power constraints. We also assume the battery is not required to achieve a final state-of-charge, but the model can easily include this constraint without loss of generality as long as the final state-of-charge is feasible. This model provides a baseline for comparing decentralized approaches, and could be extended to other DERs such as electric vehicles in the context where the generation is zero marginal-cost and energy constraints are relevant. In the case where the DERs introduce additional costs or utility functions, such as fuel costs or non-concave utility functions, the theoretical results may not hold.

The key theoretical insight we provide is that in the presence of energy storage, the dispatch cannot be controlled by price alone. Specifically, we show that if individuals act independently to maximize their utility in the presence of an optimal price, there is no guarantee that their corresponding target power injections will be feasible and satisfy power balance. This highlights that ensuring feasibility is an important requirement of decentralized mechanisms. We describe why this is not trivial in the presence of storage, and also derive equations describing the optimal power and price trajectories.

\subsection{Utility maximization model}

The model \eqref{eq:centralized} is similar in structure to a standard discrete-time, centralized energy management system. The central constraint is matching supply and demand on the time scale of hours, while we assume that droop-like control of power converters is necessary and sufficient to adjust any power imbalance in the short-term.\footnote{This technology has been extensively studied, and while important questions remain for large system stability with high penetrations of converter-interfaced generation, a variety of techniques have been validated for microgrids, see e.g. \cite{rocabert2012control}.} We include operational constraints on energy storage, but not the network constraints,\footnote{The model can be extended to include linearized power flow and line loading constraints, which would add some complexity without affecting the main results; however, full AC power flow equations would destroy the constraint linearity (and convexity) that the analysis relies on.} and assume strictly concave utility functions $U_{n,t}$. 
\begin{subequations}
\label{eq:centralized}
\begin{align}
    &\min_{\boldsymbol{p}, \boldsymbol{d}, \boldsymbol{s}} ~ - \sum_{t\in \mathcal{T}} \sum_{n\in \mathcal{C}} U_{n,t} (d_{n,t}) \\
    &\text{ s.t. } ~~ \pi_t: \sum_{n\in \mathcal{C}} d_{n,t} = \sum_{i \in \mathcal{B}} p_{i,t}^b + \sum_{g\in \mathcal{G}} p_{g,t}^s, \, \forall t\in \mathcal{T} \label{eq:power_balance} \\
    &\hspace{0.9cm} \lambda_{g,t}^{s}: 0 \le p_{g,t}^{s} \le 
\bar{P}_{g,t}^s, \, \forall g \in \mathcal{G}, \forall t \in \mathcal{T} \label{eq:pv_lims}\\
&\hspace{0.9cm} {\lambda_{n,t}^{d,-}:-d_{n,t} \le 0, \, \forall n \in \mathcal{C}, \forall t \in \mathcal{T} \label{eq:demand_pos}}\\
&\hspace{0.9cm} \lambda_{i,t}^{b}:-\bar{P}_{i,t}^b \le p_{i,t}^b \le \bar{P}_{i,t}^b, \, \forall i \in \mathcal{B}, \forall t \in \mathcal{T} \label{eq:batt_rate}\\
&\hspace{0.9cm} \lambda_{i,t}^{c}:  0 \le s_{i,t}\le 
\bar{S}_{i,t}, \, \forall {i} \in \mathcal{B}, \forall t \in \mathcal{T} \label{eq:batt_lims} \\
&\hspace{1.9cm} s_{i,t} = s_{i,t-1} - p_{i,t}^b \Delta T, \, \forall {i}\in \mathcal{B}, \forall t \in \mathcal{T}. \label{eq:SOC}
\end{align}
\end{subequations}

This allows battery constraints to be time-varying but typically $\bar{P}^b$ and $\bar{S}$ are static. The dual variables of the respective constraints are indicated before the colon. For compactness, we use a single variable to represent the difference in upper and lower bound duals, $\lambda := \lambda^+ - \lambda^-$. The initial state of charge $s_{i,0}$ is a parameter. We eliminate the constraint \eqref{eq:SOC} and decision variables $s_{i,t}$ by solving for it as $s_{i,t} = s_{i,0}-\Delta T\sum_{\tau\leq t}p^b_{i,\tau}$ and substituting this into \eqref{eq:batt_lims}.

\subsection{Theoretical analysis}
Firstly, note that all constraints in \eqref{eq:centralized} are affine, thereby satisfying the linearity constraint qualification (LCQ). This implies that for a locally optimal primal solution, there exists a set of dual variables satisfying the Karush-Kuhn-Tucker (KKT) conditions. Secondly, as all $U_{n,t}$ are concave, the problem is convex. Any point satisfying the KKT conditions is thus globally optimal and strong duality holds.

\begin{remark}
(Dual decomposition into private decisions): The Lagrangian dual of the centralized problem \eqref{eq:centralized} is separable and equivalent to the sums of Lagrangian duals for constrained individual welfare maximization for a price equal to $\pi_t$. This allows interpretation of $\pi_t$ as the electricity price. Assuming the  utility functions are concave, the Lagrangian dual problem gives the optimal price and total welfare.
\end{remark}
The Lagrangian of \eqref{eq:centralized} can be written as:
\begin{align*}
   & \mathcal{L}(d,p^s,p^b,\pi,\lambda) = \sum_{t\in \mathcal{T}}
   \sum_{n \in \mathcal{C}} \Big(\!\!-U_{n,t}(d_{n,t}) \!+ \!(\pi_t-{\lambda_{n,t}^{d,-}})\, d_{n,t}\!\Big) \nonumber\\
   &+ \sum_{g \in \mathcal{G}} \Big((\lambda_{g,t}^s-\pi_t) p_{g,t}^s + \lambda_{g,t}^{s,+} \bar{P}^s_{g,t} \Big)
   + \sum_{i \in \mathcal{B}} \Big((\lambda_{i,t}^b-\pi_t) p_{i,t}^b \nonumber\\
   &+ \lambda^c_{i,t}\Big(s_{i,0} -\Delta T \sum_{\tau \leq t} p_{i,\tau}^b\Big) - (\lambda^{b,+}_{i,t} + \lambda^{b,-}_{i,t})\bar{P}^b_{i,t} - \lambda^{c,+}_{i,t}\bar{S}_{i,t} \Big). \label{eq:centralized_lagrangian}
\end{align*}
We define individual utility/profit-maximization problems for each of the consumers, PV, and battery operators for an electricity price as in \eqref{eq:W_nt}-\eqref{eq:W_it}.
\begin{align}
    W_n(\pi) &:= \min_{d_{n} \ge 0} \sum_t -U_{n,t}(d_{n,t}) + \pi_t\, d_{n,t}, \label{eq:W_nt} \\
    W_g(\pi) &:= \min_{p_{g}^s} \sum_t -\pi_t \,p_{g,t}^s \quad \text{s.t.} \quad \eqref{eq:pv_lims}, \label{eq:W_gt}\\
    W_i(\pi) &:= \min_{p_{i}^b} \sum_t -\pi_t\, p_{i,t}^b \quad \text{s.t.} \quad \eqref{eq:batt_rate}-\eqref{eq:SOC}. \label{eq:W_it}
\end{align}
Denoting their Lagrangians by $\mathcal{L}_n$, $\mathcal{L}_g$, $\mathcal{L}_i$, one can show that
\begin{align}
    \mathcal{L}(d,p^s,p^b,\pi,\lambda) =& \sum_{n \in \mathcal{C}} \mathcal{L}_n(d_n,\pi) + \sum_{i \in \mathcal{G}} \mathcal{L}_g(p_g^s,\lambda_g^s,\pi) \nonumber \\
    &+ \sum_{i \in \mathcal{B}} \mathcal{L}_i(p_i^b,\lambda_i^b,\lambda_i^c,\pi).
\end{align}
As $W_g$ and $W_i$ are linear programs, strong duality holds for these subproblems, and the Lagrangian dual problem is
\begin{align}
    &\max_{\pi,\lambda} \inf_{d,p^s,p^b} \mathcal{L}(d,p^s,p^b,\pi,\lambda)\nonumber\\ = &\max_\pi \sum_{n \in \mathcal{C}} W_n(\pi) + \sum_{g \in \mathcal{G}} W_g(\pi) + \sum_{i \in \mathcal{B}} W_i(\pi). \label{eq:dual_prob}
\end{align}
By strong duality \eqref{eq:dual_prob} gives the optimal objective value with its maximizer $\pi^\star$ equal to the optimal price. However, as we establish later, the optimal $p_i^b$ for \eqref{eq:W_it} is not necessarily unique, meaning that broadcasting an optimal price to individual agents does not necessarily satisfy constraint \eqref{eq:power_balance} and clear the market; i.e., primal feasibility is not guaranteed.

\begin{remark}
For all $t \in \mathcal{T}$, the following relations hold true at optimum and characterize the optimal price
\end{remark}
\begin{subequations}
\begin{align}
    \pi^\star_t& = {\partial U_{n,t}(d^\star_{n,t})}/{\partial d_{n,t}}+\lambda_{n, t}^{\star,d,-}, \, \forall n \in \mathcal{C} \label{eq:statio_d}\\
&=\lambda_{i,t}^{\star,b}\!-\! \Delta T\sum_{\tau\geq t}\!\lambda_{i,\tau}^{\star,c}, \, \forall i \in \mathcal{B} \label{eq:statio_batt}\\
&=\lambda_{g,t}^{\star,s}, \, \forall g \in \mathcal{G}.
\end{align}
\end{subequations}

Each of the equalities follow from the stationarity conditions of \eqref{eq:centralized}. We interpret the dual variable $\pi^\star_t$ as the price by \textbf{Remark~1} and note from \eqref{eq:statio_batt} that it depends on the cumulative future shadow prices of the storage capacity constraint. Eq.~\eqref{eq:statio_d} requires $U_{n,t}$ to be differentiable for equality but can be replaced by the subdifferential of $U_{n,t}$ otherwise.

\begin{remark}
If at time $t$, a utility function for at least one customer is differentiable and strictly increasing on $\mathbb{R}^+$, then at optimum, the price is strictly positive and solar production is at its maximum.
\end{remark}
This follows from \textbf{Remark~2} and the properties of strictly increasing functions:
\begin{equation*}
\exists n \ni {\partial U_{n,t}(d^\star_{n,t})}/{\partial d_{n,t}} > 0 \; \forall d_{n,t}\Rightarrow \pi^\star_t > 0 \Rightarrow \lambda_{g,t}^{\star,s} > 0.
\end{equation*}
By complementary slackness, $\lambda_{g,t}^{\star,s} > 0 \Rightarrow p_t^{\star,s}=\bar{P}_t^s$. This is intuitive as it is better to supply any benefiting consumer than curtailing available solar. This also implies that solar generation can be removed as a decision variable and set to the available resource in this case.

\begin{remark}
The optimal price evolves as
\begin{equation}
\label{eq:rem4}
\pi^\star_{t+1}-\pi^\star_t = \lambda_{i,t+1}^{\star,b}-\lambda_{i,t}^{\star,b}+\Delta T\lambda_{i,t}^{\star,c}.
\end{equation}
\end{remark}
This follows from {\textbf{Remark~2}} by expanding the expression $\pi^\star_{t+1}\!-\!\pi^\star_t$. This captures the price trajectory, from which price volatility can be analyzed. Note that both $\lambda_{i,t}^{\star,b}$, $\lambda_{i,t}^{\star,c}$ can be less than 0. We will use \eqref{eq:rem4} for our analysis in {\bf Remark~5}.

\begin{remark}
(Non-uniqueness of decentralized battery dispatch): There are non-trivial optimal prices $\pi^\star$ such that the optimal individual battery dispatch $W_i(\pi^\star)$ is not-unique.
\end{remark}
This can be observed in a simple example. Suppose $T\!=\!5$, $\Delta T\!=\!1$, $\bar{P}^{b}_{i,t}\!\equiv\!3$, $\bar{S}_{i,t}\!=\!10$, $s_{0,i}\!=\!5$, and $\pi^\star\!=\![1,1,2,3,1]$. One can verify that $p^b_i\!=\![-0.5, -0.5, 3, 3, 0]$, $p^b_i\!=\![0, -1, 3, 3, 0]$, and $p^b_i\!=\![-1, -1, 3, 3, 1]$ are all optimal solutions with a net cost of $-14$. Here, equal prices imply there is no change in cost to shift energy from one period to another and the constraints allow this shift. More formally, if an optimal solution is not on any of the constraint boundaries \eqref{eq:batt_rate}-\eqref{eq:batt_lims} at $t$ and $t\!+\!1$, then it will not be unique because not being on the boundary implies 1) that $\lambda^{\star,b}_{i,t}$, $\lambda^{\star,b}_{i,t+1}$, $\lambda^{\star,c}_{i,t}$, and $\lambda^{\star,c}_{i,t+1}$ are all $0$, so $\pi^\star_{t+1}\!=\!\pi^\star_t$ by \eqref{eq:rem4}, therefore $p'^{,b}_{i,t}\!=\!p^{\star,b}_{i,t}\!+\!\varepsilon$ and $p'^{,b}_{i,t+1}\!=\!p^{\star,b}_{i,t+1}-\varepsilon$ have equivalent net cost $\forall \varepsilon$ without affecting the solution at other times; and 2) that this perturbation is feasible for sufficiently small $|\varepsilon|\!>\!0$. Note that this condition is overly restrictive and not necessary for non-uniqueness; in particular energy may be shifted between non-consecutive time periods, and only particular combinations of constraints between the time periods need to be non-binding rather than all constraints. Equal prices between time periods may indicate non-uniqueness, but the optimal solution may still be unique if the constraints do not allow a perturbation to remain feasible. In Appendix~\ref{sec:detailed_battery_model}, we show that this applies with a variation in the optimal price profile even when battery inefficiencies and self-discharge are considered.

The consequence of \textbf{Remark~5} is that in general, an optimal price is not sufficient to yield individual battery dispatch decisions with the optimal quantity, meaning that a system operator cannot control dispatch outcomes solely by broadcasting a price signal or adequately forecast the decentralized response to price. Even when all utility functions are strictly concave so that the solution to the centralized problem \textit{is} unique and corresponds to an optimal price $\pi^\star$, there are (likely common) conditions whereby an individual battery operator's decision in response to $\pi^\star$ does not satisfy the constraint \eqref{eq:power_balance}. Intuitively, if the price is constant between two successive periods, a battery operator would be indifferent to selling more energy in one period versus the next, so their dispatch is not unique and there is no guarantee that the dispatch will meet demand.

The previous observation implies that extending standard centralized market mechanisms to systems with energy storage faces limitations. If a centralized energy market limits entities with storage to submitting a single curve of price and quantity for each time-period, it is likely to result in suboptimal outcomes to the utility maximization problem and even in infeasibility. Although not shown here, we expect this result extends to load that can be shifted without cost, and to storage models with constant charge or discharge inefficiences. P2P approaches where agents explicitly agree on quantity are a potential opportunity for addressing this challenge.
\setcounter{theorem}{0}

\begin{figure}[t]
    \centering
\tikzset{every picture/.style={scale=1.0}}%
      \input{centralized_plots.tex}
  \caption{Optimal profiles for the centralized approach with 10 agents, $\bar{P}^b_i=10$, and demand elasticity $\in [-3, -2]$ for all agents. Note: The plots for $\bar{S}_\text{tot}$=$15$, $300$ kWh overlap in (b), (c), (d).}
        \label{figure:Centralized}
\end{figure}
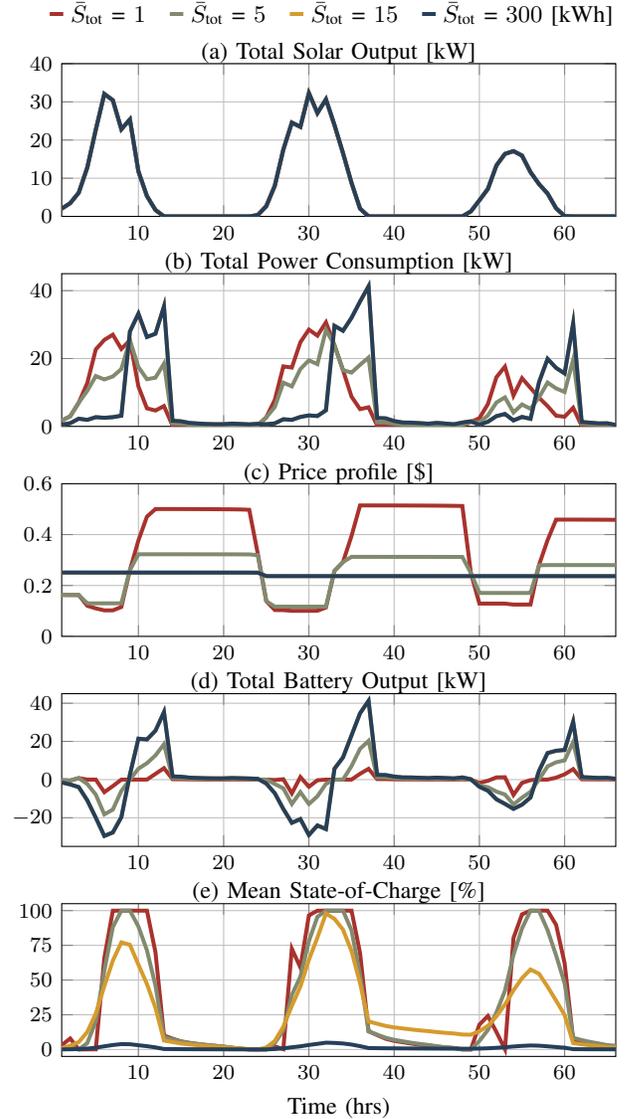


\subsection{Example optimal trajectories and the effects of storage}\label{subsec:centralized_simulation}
To show how the PV profile and storage capacity affect the optimal trajectories of~\eqref{eq:centralized}, we simulate scenarios with total storage capacity varying $\bar{S}_\text{tot}\in \{1,5,15,300\}$ kWh and distributed evenly to batteries collocated with consumers. We sample hourly load and PV profiles from the $2017$ Pecan Street data set \cite{pecan2020dataport} over a random $66$ hour interval and construct example utility functions by assuming a quasi-constant price-elasticity demand function and centering it at the observed load with a constructed time-varying price profile\footnote{$0.10$/kWh between $21$:$00$-$11$:$00$, $0.15$/kWh between $11$:$00$-$16$:$00$, and $0.30$/kWh between $16$:$00$-$21$:$00$.} (see Appendix~\ref{sec:utility_function}).We model $10$ consumers and randomly select elasticities $\in [-3,\,-2]$.

Figure~\ref{figure:Centralized} shows the optimal trajectories for each storage scenario. The solar output is identical across scenarios (a). As the storage capacity is increased, the consumption shifts to evening peaks from daytime peaks coincident with solar (b). Increasing storage reduces the swing between high and low price periods (c). The instances when the price changes, correspond to when the battery constraints are binding, as predicted by \eqref{eq:rem4}, which also explains how higher storage capacities lead to a flat price by reducing $\lambda_{i,t}^{\star,c}$ to a negligible value. A flatter price means the conditions of \textbf{Remark 5} are less likely to be met, highlighting the increasing need to coordinate battery dispatch as capacity increases. In contrast, a smaller capacity induces cyclical price fluctuations through peak-to-peak cycling. This also illustrates how the marginal value of storage in arbitraging high and low price periods depends on the existing capacity. These phenomena are explained analytically by the model; extending the model to derive optimal investment and planning decisions is a promising area for future work.

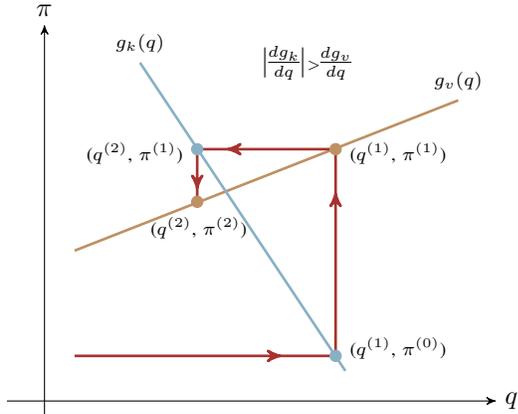
\begin{figure}[t]
\tikzset{every picture/.style={scale=1.0}}%
\centering
\input{cobweb.tex}
\caption{The standard convergent cobweb model.}
\label{fig:cobweb}
\end{figure}

\section{Peer-to-peer Negotiation}

In this section we analyze how a decentralized, peer-to-peer energy market can arrive at a near-optimal dispatch solution using an intuitive negotiation approach. We model a process of exchanging price and quantity offers after the classic ``cobweb'' model of dynamic markets \cite{kaldor1934classificatory} and observe that classical results show the process can diverge. We therefore, consider an additional \textit{dynamic step-limiting constraint} on the process to ensure convergence, which could be thought of as a behavioral tendency of agents or an explicit rule to be imposed by a bidding platform. We assume agents are matched \textit{a priori} and that offers are synchronized so as to simplify the analysis and presentation, but posit that the process can be generalized to capture more informal interaction between agents.

As a starting point, consider an interaction between two agents who are ``prosumers'' with private solar and storage systems and who individually derive private value from energy use. Most likely, there exists a trade that makes both agents better off. An intuitive way for the agents to find such a trade is for one to start by proposing a quantity (either positive or negative) and for the other to respond with a price. The first agent would likely reassess the quantity they would seek at that price, propose a new quantity, and so on. This iterative process is described by the cobweb model illustrated in Fig.~\ref{fig:cobweb}. The equilibrium is the intersection of supply and demand curves arising from the utility functions. This is the optimum of the utility maximization model but the process converges to this point if and only if the magnitude of the slope of the demand curve exceeds that of the supply curve at the equilibrium \cite{kaldor1934classificatory}.

\begin{figure}[t]
    \centering
\tikzset{every picture/.style={scale=1.0}}%
      \input{probnoncon.tex}
  \caption{Convergent trajectory under the dynamic step-limiting constraint where the standard cobweb model would diverge.}
        \label{fig:noncon}
\end{figure}
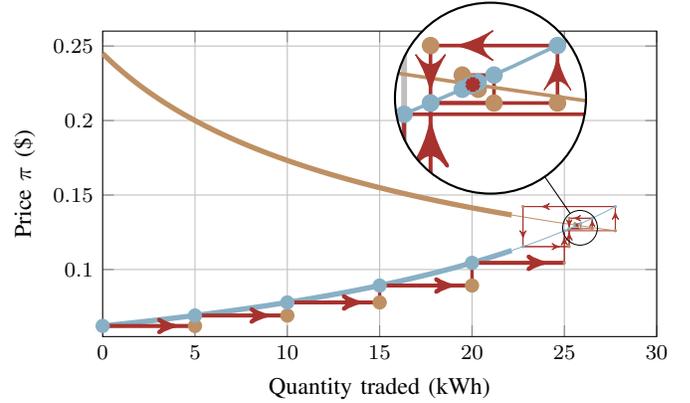

We modify the cobweb model to ensure convergence even when this condition is not met by including a step-limiting constraint, illustrated in Fig.~\ref{fig:noncon}. This constraint assumes (or enforces) that agents will not adjust their quantity offers by more than some threshold each iteration, and that this threshold shrinks if the quantity is ``oscillating''. We generalize to consider multiple agents proposing quantities (called $q$-agents) to agents who respond with price (called $\pi$-agents). The agents exchange vectors of quantity and price for each period over a finite time horizon. To simplify the analysis, we assume a single $\pi$-agent interacts with multiple $q$-agents. In practice, there would likely be multiple $\pi$-agents, and $q$-agents would select one or more $\pi$-agents to negotiate with, based on their expectation of the outcome of the negotiation, but this matching problem is beyond the scope of this paper.

We present formal decision models for the $q$-agents and the $\pi$-agents, and define an iterative process that guarantees physically feasible and at least weakly Pareto-optimal outcomes (i.e., no agents are worse off). We prove theoretically that the process converges to within a tolerance of the centralized solution for the $2$-agent, single time step case, and demonstrate convergence using simulations for the general case in the next section. These results show that an informal, decentralized, peer-to-peer negotiation process is capable of approximating the centralized welfare maximization problem, and offers a specific approach that could be implemented on a software platform and evaluated in practice. In addition, and contrary to the centralized approach, this negotiation process does not require full information exchange between agents, since private utility functions are hidden and only trading of quantities and prices are required for the execution of the algorithm.

We denote the set of $\pi$-agents with $\mathcal{V}$ and $q$-agents with $\mathcal{U}$, such that $\mathcal{C} = \mathcal{U}\cup \mathcal{V}$ and $\mathcal{U} \cap \mathcal{V} = \varnothing$. We index the $q$-agents by $k \in \mathcal{U}$ and the single $\pi$-agent as $v$, $\mathcal{V} = \{v\}$. The $q$-agents may exit the process early, which we track by partitioning $\mathcal{U}$ into exited agents $\mathcal{X}$ and negotiating agents $\mathcal{Y}$, and updating these dynamically. All constraints are implicitly defined $\forall t \in \mathcal{T}$.

\subsubsection{Optimization problem for the $\pi$-agent}

The $\pi$-agent receives a set of requested quantities $\boldsymbol{q}\!=\!\{\boldsymbol{q}_{k}\}$ from each $k \!\in\! \mathcal{Y}$ (positive means $k$ receives energy), with $\boldsymbol{q}_k \!=\! \{\boldsymbol{q}_{k,t}\}$, which may not be feasible. The $\pi$-agent first projects $\boldsymbol{q}$ to a feasible $\boldsymbol{q}'$ by keeping reference to a quantity $\boldsymbol{\hat{q}}$ known to be feasible to all agents; $\boldsymbol{q}'$ is restricted to lie on the line connecting $\boldsymbol{q}$ and $\boldsymbol{\hat{q}}$ defined by \eqref{eq:projection_line}, where $\beta\!=\!0$ yields the requested $\boldsymbol{q}$ and $\beta\!=\!1$ the known feasible $\boldsymbol{\hat{q}}$. Thus, minimizing $\beta \!\geq\! 0$ selects the closest point to $\boldsymbol{q}$ satisfying the constraints:
\begin{subequations}
 \begin{align}
   &\min_{\boldsymbol{d}_v, \boldsymbol{p}_v^b, \boldsymbol{s}_v, \boldsymbol{q}', \beta}~ \beta \\
    &\text{ s.t. }~~ d_{v,t} + \sum_{k \in \mathcal{Y}}q'_{k,t} + \sum_{k \in \mathcal{X}} q_{k,t} = p_{v,t}^s + p_{v,t}^b \\
    &\hspace{1.2cm} 0 \le \beta \le 1 \\
    &\hspace{1.2cm} q'_{k,t} = \beta\,\hat{q}_{k,t} + (1-\beta)\,q_{k,t} \label{eq:projection_line} \\
    &\hspace{1.2cm}\text{and constraints } \eqref{eq:demand_pos}-\eqref{eq:SOC}.
\end{align}
\end{subequations}
We maintain that $\boldsymbol{\hat{q}}$ is feasible for all agents. Before any agents exit, $\mathcal{X} \!=\! \varnothing$ and $\boldsymbol{\hat{q}}\!=\!0$ is feasible, so we initialize with $\boldsymbol{\hat{q}}\!=\!0$ and update $\boldsymbol{\hat{q}}$ as agents exit at feasible points. As shown below, $\boldsymbol{q}$ is necessarily feasible for each $q$-agent, and their constraints are convex, so any point connecting two feasible points is feasible, and in particular $\boldsymbol{q}'$. It is also possible to include additional constraints in this optimization problem, for example, to ensure power flow feasibility. Extending this model to include a network feasibility validation as constraints over $q_{k,t}'$ is a promising direction for future work.

Next, the $\pi$-agent solves their utility maximization problem to obtain $\boldsymbol{\pi}\!=\!\{\pi_t\}$ and their utility from these proposed trades. A key assumption is that they set $\boldsymbol{\pi}$ at their marginal utility; i.e., they bid according to a competitive market strategy and cannot exercise market power. This is likely to hold in practice if there are sufficiently many $\pi$-agents the $q$-agents can access; however, we recommend a more careful analysis of market power in the scope of a ``many-to-many'' extension to this work. The maximization problem is:
\begin{subequations}
 \begin{align}
   &\min_{\boldsymbol{d}_v, \boldsymbol{p}_v^b, \boldsymbol{s}_v} ~ - \sum_{t\in\mathcal{T}}  U_{v,t}(d_{v,t})\\
   &\text{ s.t. }~~ ~\pi_t: d_{v,t} + \sum_{k \in \mathcal{Y}}q'_{k,t} + \sum_{k \in \mathcal{X}} q_{k,t} = p_{v,t}^s + p_{v,t}^b \label{eq:dual_pi}\\
   &\hspace{1.2cm}\text{and constraints } \eqref{eq:demand_pos}-\eqref{eq:SOC}.
\end{align}
\end{subequations}
As in the centralized model, the price is given directly by the stationarity condition with $\lambda_v^{d,-}\!=\!0$:
$$\pi_t = {\partial U_{v,t}({d^\star_{v,t}})}/{\partial d_{v,t}}.$$
Lastly, the $\pi$-agent checks whether its utility from this potential trade is at least as high as its optimal utility from no trade (specifically solving the same problem with $\boldsymbol{q}'\!=\!0$), and sets a binary variable $\alpha_v$ true if so, and false otherwise. This $\alpha_v$ signals whether $v$ would prefer $\boldsymbol{q}'$ to no trade. We denote the entire decision as $\mathcal{P}_v^\pi\!: \, (\boldsymbol{q},\boldsymbol{\hat{q}}) \!\mapsto\! (\boldsymbol{q}',\boldsymbol{\pi},\alpha_v)$.

\subsubsection{Optimization problem for $q$-agents}

The $k$-th $q$-agent makes the decision $\mathcal{P}_k^q\!: \, (\boldsymbol{\pi},\boldsymbol{q}'_k, \boldsymbol{\delta}_k) \!\mapsto\! (\boldsymbol{q}_k,\alpha_k,\eta_k)$, where $\alpha_k$ carries the analogous meaning to $\alpha_v$, $\eta_k$ signals whether they are ``satisfied'', $\boldsymbol{q}'_k$ is the subset of $\boldsymbol{q}'$ for $k$, and $\boldsymbol{\delta}_k$ is the step-limiting constraint restricting the $q$-agent to select something close to the offer $\boldsymbol{q}'$. The decision is: 
\begin{subequations}
 \begin{align}
    &\min_{\boldsymbol{d}_k, \boldsymbol{q}_k, \boldsymbol{p}_k^b, \boldsymbol{s}_k}~   \sum_{t\in\mathcal{T}}  -U_{k,t}(d_{k,t}) + \pi_t q_{k,t}  &\\
    &\text{ s.t. }~~~  d_{k,t} - p_{k,t}^s - p_{k,t}^b - q_{k,t} = 0
    \label{eq:bal_eq2}\\
    &\hspace{1.0cm} \left|q_{k,t} - q_{k,t}'\right| \le \delta_{k,t}
    \label{eq:bat_step2}\\
    &\hspace{1.1cm}\text{and constraints } \eqref{eq:demand_pos}-\eqref{eq:SOC}.
\end{align}
\end{subequations}
Agent $k$ requests to finalize the trade and exit if their (not necessarily unique) optimal $\boldsymbol{q}_k$ is close enough to the offer $\boldsymbol{q}'_k$, where the distance is determined by a small $\varepsilon$:
\begin{align}
    \eta_k = \left\{ \begin{array}{cl}
        \text{True} & \text{if } |q_{k,t} - q'_{k,t}| \le \gamma\, \varepsilon \\
                    \text{False}  & \text{otherwise}.
    \end{array} \right.
\end{align}
The exit condition includes the constant $\gamma \!\in\! (0,1)$ to simplify the statement of Theorem~6, but could be modified with an update to the bound in the theorem. An alternative criterion based on whether the utilities from these offers are close enough could also be used but would affect the bound.

\subsubsection{Iterative Algorithm}
\begin{algorithm}[h]
\SetAlgoLined
\KwResult{Energy trades ($\boldsymbol{\pi}^\star_k, \boldsymbol{q}^\star_k$) for each agent $k \!\in\! \mathcal{C}$.}
 \textbf{Initialization:}
  Define the $\pi$-agent $v \!\in\! \mathcal{C}$ and the parameters $\gamma \!\in\! (0,1)$, $\varepsilon \!>\! 0$, initial step-limit $\delta^{(0)} \!>\! \gamma \varepsilon$\, and max iterations $M$\;
  Set $i \!\leftarrow\! 1$, $(\boldsymbol{q}^{(1)},\boldsymbol{\hat{q}})\!\leftarrow\!(0,0)$, $\{\delta^{(1)}_{k,t}\}\!\leftarrow\!\delta^{(0)}$, $\mathcal{X}\!\leftarrow\!\{0\}$, and $\mathcal{Y}\!\leftarrow\!\mathcal{C} \setminus{\{v\}}$ \;
  \While{$\mathcal{Y} \!\neq\! \varnothing$ \textbf{and} $i \!\leq\! M $}{
  $(\boldsymbol{q}^{\prime,(i)},\boldsymbol{\pi}^{(i)},\alpha_v)\!\leftarrow\!\mathcal{P}_v^\pi(\boldsymbol{q}^{(i)},\boldsymbol{\hat{q}})$\;
  \For{$k \!\in\! \mathcal{Y}$}{
    $(\boldsymbol{q}_k^{(i+1)},\alpha_k,\eta_k) \!\leftarrow\! \mathcal{P}_k^q(\boldsymbol{\pi}^{(i)},\boldsymbol{q}_k^{\prime,(i)},\boldsymbol{\delta}_k^{(i)})$\;
    $\boldsymbol{o}_k^{(i)}\!\leftarrow\!f^o(\boldsymbol{q}_k^{(i+1)},\boldsymbol{q}_k^{(i)},\boldsymbol{q}_k^{(i-1)})$\;
    \leIf{$\eta_k$}{$\boldsymbol{\delta}_k^{(i+1)} \!\leftarrow\! \boldsymbol{\delta}_k^{(i)}$}{$\boldsymbol{\delta}_k^{(i+1)} \!\leftarrow\! f^\delta (\boldsymbol{\delta}_k^{(i)},\boldsymbol{o}_k^{(i)})$}
    }
    \If{$\alpha_j \; \forall j \!\in\! \mathcal{Y} \!\cup\! \{v\}$}{
        $\boldsymbol{\hat{q}} \!\leftarrow\! \boldsymbol{q}^{\prime,(i)}$\;
        \For{$k \!\in\! \mathcal{Y}$ \textbf{where} $\eta_k$}{
            $\mathcal{Y} \!\leftarrow\! \mathcal{Y} \!\setminus\! \{k\}$, $\mathcal{X} \!\leftarrow\! \mathcal{X} \!\cup\! \{k\}$, $(\boldsymbol{\pi}^\star_k, \boldsymbol{q}^\star_k) \!\leftarrow\! (\boldsymbol{\pi}^{(i)},\boldsymbol{q}_k^{\prime,(i)})$
        }
    }
$i \!\leftarrow\! i + 1 $
  }
  \caption{Bounded cobweb iteration for a single $\pi$-agent and multiple $q$-agents.}
  \label{alg:4}
\end{algorithm}

The negotiation algorithm is presented in Algorithm \ref{alg:4}. At each iteration, $q$-agents submit their energy quantity requests to the $\pi$-agent based on the last price and quantity offered by the $\pi$-agent. The $q$-agents are allowed to exit only when all agents have declared $(\boldsymbol{\pi}, \boldsymbol{q}')$ preferable to no trade through $\alpha$ (i.e., $\alpha_k$ \!=\! True $\forall \, k \!\in\! \mathcal{U}$), guaranteeing that trades are at least weak-Pareto improvements. Importantly, the step-limit $\boldsymbol{\delta}$ is shrunk by $\gamma \!\in\! (0,1)$ if the quantity is ``oscillating'' (see Fig.~\ref{fig:cobweb}), defined by the binary state $\boldsymbol{o}^{(i)}$ as the quantity not monotonically increasing or decreasing over 3 iterations, with $o^{(1)} \!=\! 1$ and update maps $f^o$ and $f^\delta$:
\begin{align*}
    &f^o:o_{k,t}^{(i)}\!= \!\neg (q_{k,t}^{(i+1)}\! > \!q_{k,t}^{(i)} \!>\!q_{k,t}^{(i-1)} \vee q_{k,t}^{(i+1)}\! <\! q_{k,t}^{(i)}\! < \!q_{k,t}^{(i-1)}), \\
    &f^\delta:\delta^{(i+1)}_{k,t} = (1-o_{k,t}^{(i)})\,\delta^{(i)}_{k,t} + o_{k,t}^{(i)}\,\gamma\delta^{(i)}_{k,t}. \label{eq:delta_update}
\end{align*}
This shrinking step-limit prevents the divergent case of the cobweb model \cite{kaldor1934classificatory}.

\subsubsection{Optimality of the two-agent, single time step case}
In this subsection we prove that Algorithm \ref{alg:4} converges within an $\varepsilon$ tolerance in finite iterations to the socially optimal quantity in the case of only two agents with single time horizon. We ignore storage in this case, as it can equivalently be treated as solar production for $\mathcal{T} \!=\! \{1\}$, and drop the time index $t$ for brevity. We assume the the solar production is greater than zero for at least one agent, and that each agent's marginal utility of consumption $\partial U_n(d_n)/\partial d_n$ is strictly monotonically decreasing on $[0,\infty)$ and decreasing asymptotically to zero.

Note that $q \!=\! d_k \!-\! p_k^s \!=\! -d_v \!+\! p_v^s$, and the \textit{unconstrained} demand and supply curves are defined as $g_k\!\equiv\!\partial U_k(q)/\partial d_k$ and $g_v\!\equiv\!\partial U_v(q)/\partial d_v$. Thus, $g_k$ is monotonically decreasing and $g_v$ is monotonically increasing. Without the step-limiting constraint \eqref{eq:bat_step2}, the problem $\mathcal{P}^q_k$ for the $q$-agent has a closed form solution:
\begin{equation}
    q^\dagger = g_k^{-1}(\min (g_k(-p_k^s),\pi)) \equiv h_k(\pi),\label{eq:A_map}
\end{equation}
where $g_k^{-1}$ denotes the inverse of $g$ with domain $(0,g_k(-p_k^s)]$. With the step-limiting constraint, the solution is
\begin{align}
    q = \left\{\begin{array}{cc}
        q^\dagger & \text{ if } |q^\dagger - q'| \le \delta  \\
        q' + \delta  & \text{ if } q^\dagger >  q' + \delta \\
        q' - \delta  & \text{ if } q^\dagger <  q' - \delta.  
    \end{array} \right. \label{eq:q_sol}
\end{align}

The projection step reduces to $q' \!=\! \min(p^s_v,q)$, and the $\pi$-agent's price is given by $\pi \!=\! g_v(q')$.

The optimal quantity of the centralized problem $q^\star$ is the unique fixed point of the iteration if $q^\star < p_v^s$ or if $q^\star \!=\! p_v^s$ and $g_k(p_v^s) \!=\! g_v(p_v^s)$. Indeed, note that $-p_k^s \!\leq\! q^\star \!\leq\! p_v^s$ by the constraints, and that $h_k(\pi^\star) \!\equiv\! q^\star$. If $q^{(i)} \!=\! q^\star$, then $q' \!=\! q^\star$ and $\pi^{(i)} \!=\! g_v(q^\star) \!=\! \pi^\star \!-\! \lambda^{\star,d,-}_s$ by \eqref{eq:statio_d}. When $q^\star \!<\! p_v^s$ or $g_k(p_v^s) \!=\! g_v(p_v^s)$, then we have $\lambda^{\star,d,-}_s \!=\! 0$ and $\pi^{(i)} \!=\! \pi^\star$, and hence $q^\dagger \!=\! q^\star$ with $q^{(i+1)} \!=\! q^\star$. Otherwise, $\lambda^{\star,d,-}_s \!>\! 0$ and $\pi^{(i)} \!<\! \pi^\star$, so $q^\dagger \!>\! q^\star$ by the strict monotonicity of $h_k$, and $q^{(i+1)} \!>\! q^\star$, so it is not a fixed point. In other words, the fixed point is the intersection of the curves $g_v$ and $g_k$, as shown in Fig. \ref{fig:cobweb}. Since both curves are strictly monotonic, this fixed point is unique. If they do not intersect on $[-p_k^s,p_v^s]$, then $q^\star$ is only a fixed point if it is $-p_k^s$.

\begin{lemma}[Movement towards equilibrium]\label{lemma:movement} At any iteration $i$, if $q^\star \!<\! p_v^s$, then $q' \!\leq\! q^\star \!\Leftrightarrow\! q^{(i+1)} \!\geq\! q'$ and $q' \!\geq\! q^\star \!\Leftrightarrow\! q^{(i+1)} \!\leq\! q'$. Moreover, $q' \!\leq\! q^\star  \!\Leftrightarrow\! q^{(i+1)} \!\geq\! q^{(i)}$ and $q' \!\geq\! q^\star \!\Leftrightarrow\! q^{(i+1)} \!\leq\! q^{(i)}$.
\end{lemma}
\textbf{Proof:} We prove this by showing the forward direction of the first set of statements $q' \!\leq\! q^\star \!\Rightarrow\! q^{(i+1)} \!\geq\! q'$ and $q' \!\geq\! q^\star \!\Rightarrow\! q^{(i+1)} \!\leq\! q'$. Each of these statements implies the converse of the other is true, establishing the reverse direction. We use the same approach for the second set of statements.

Let $\pi \!=\! g_v(q')$. If $q^\star \!<\! p_v^s$, then $g_v(q^\star) \!=\! \pi^\star$ and if $q' \!\leq\! q^\star \!\Rightarrow\! \pi \!\leq\! \pi^\star \!\Rightarrow\! h_k(\pi) \!\geq\! h_k(\pi^\star) \!\Rightarrow\! q^\dagger \!\geq\! q^\star$ because $g_v$ is monotonically increasing and $h_k$ is monotonically decreasing. Thus, $q^\dagger \!\geq\! q^\star \!\geq\! q' \!\Rightarrow\! q^{(i+1)} \!=\! \min(q^\dagger,q'+\delta) \!\Rightarrow\! q^{(i+1)} \!\geq\! q'$. By the same logic, $q' \!\geq\! q^\star \!\Rightarrow\! q^{(i+1)} \!\leq\! q'$. Moreover, $q' \!\leq\! q^\star \!<\! p^v_s \!\Rightarrow\! q' \!=\! q^{(i)}$; therefore, $q^{(i+1)} \!\geq\! q' \!\Rightarrow\! q^{(i+1)} \!\geq\! q^{(i)}$. Finally, because $q' \!=\! \min(p_v^s, q^{(i)}) \!\leq\! q^{(i)}$, and by showing that $q' \!\geq\! q^\star \!\Leftrightarrow\! q^{(i+1)} \!\leq\! q'$, it follows that $q' \!\geq\! q^\star \!\Rightarrow\! q^{(i+1)} \!\leq\! q^{(i)}$.
\hfill $\square$

\begin{lemma}[Entry to the oscillatory mode] \label{lemma:oscillatory}
If the system is not in the oscillatory mode at iteration $i$, then $\exists\, l\!>\!0$ such that if the algorithm does not terminate at iteration $s \!<\! i\!+\!l$, it will be in the oscillatory mode at $i\!+\!l$. 
\end{lemma}
\textbf{Proof:} First, consider the case $q^\star \!<\! p_v^s$. We will show the case when $q^{(i)} \!<\! q^\star$, since the other case is analogous. By Lemma \ref{lemma:movement}, $q^{(i+1)}$ moves towards the equilibrium and $\delta^{(i+1)}$ is not reduced when moving in the same direction. Thus, for some $j \!>\! i$, $q^{(j)} \!\ge\! q^\star$ (with $q^{(j-1)} \!\le\! q^\star$); hence, by Lemma \ref{lemma:movement}, $q^{(j+1)} \!\le\! q^{(j)}$ and we enter the oscillatory mode at $j+1 \!=\! i\! + \!l$. 

Second, we consider the case of $q^\star \!=\! p_v^s$. Observe that when eventually $q^{(j)} \!\ge\! p_v^s$, it will be projected back to $q^{\prime} \!=\! \min(q^{(j)}, p_v^s) \!=\! p_v^s$ to ensure feasibility for the $\pi$-agent. Then, since there is no intersection of marginal utility curves in the interior, it implies that $g_k(q') \!\ge\! g_v(q') \!=\! \pi^{(j)}$, and hence, the $q$-agent again requests $q^{(j+1)} \!\ge\! p_v^s$, that gets projected back to $q^{\prime} \!=\! p_v^s$. Thus, since repeated values of $q$ are received, it will enter the oscillatory mode (and eventually converge to $p_v^s$).
\hfill $\square$

\begin{lemma}[Boundedness of distance from the equilibrium] \label{lemma:bounded}
Assume $q^\star \!<\! p_v^s$ and suppose the system is in the oscillatory mode at iteration $i$. Then, $|q'\!-\!q^\star| \!<\! \gamma^{-1}\delta^{(i)}$.
\end{lemma}
\textbf{Proof:} Let $q'^{,(i-1)}$ denote the offer from the $\pi$-agent at the previous iteration. We first prove the case when $q' \!\geq\! q^\star$ by contradiction. To this end, assume $q' \!-\! q^\star \!\geq\! \gamma^{-1}\delta^{(i)}$. This implies $q^{(i)} \!-\! q^\star \!\geq\! \gamma^{-1}\delta^{(i)}$ because $q^{(i)} \!\geq\! q'$ by $q' \!=\! \min(q^{(i)},p_v^s)$. This in turn implies $q'^{,(i-1)} \!\geq\! q^\star$ by the step-limiting constraint \eqref{eq:bat_step2} at the previous iteration (observe $\delta^{(i-1)} \!\leq\! \gamma^{-1}\delta^{(i)})$. Then, $q'^{,(i-1)} \!\geq\! q^\star \!\Rightarrow\! q^{(i)} \!\leq\! q'^{,(i-1)} \!\leq\! q^{(i-1)}$ by Lemma \ref{lemma:movement}. For the system to be oscillating with $q^{(i)} \!\leq\! q^{(i-1)}$, either $q^{(i+1)} \!>\! q^{(i)}$ (which contradicts $q' \!\geq\! q^\star$ by Lemma \ref{lemma:movement}), or we have equality at $q^{(i-1)} \!=\! q^{(i)}$ or $q^{(i)} \!=\! q^{(i+1)}$ (which implies $q^\star \!=\! q'$, by the unique fixed point, contradicting $q'\!-\! q^\star \!\geq\! \gamma^{-1}\delta^{(i)}$).

We show the second case, $q' \!<\! q^\star$, directly. We have $q^{(i+1)} \!>\! q'$ by Lemma \ref{lemma:movement} and $q^{(i)} \!=\! q'$ because $q^\star \!\leq\! p_v^s$. Thus, the oscillating mode implies $q^{(i)} \!\leq\! q^{(i-1)}$. It holds that $q'^{,(i-1)} \!\geq\! q^\star$: if $q^{(i-1)} \!>\! p_v^s$, then $q'^{,(i-1)} \!=\! p_v^s \!\geq\! q^\star$. Alternatively, if $q^{(i-1)} \!\leq\! p_v^s$, then $q'^{,(i-1)} \!=\! q^{(i-1)} \!\geq\! q^{(i)}$ and $q^{(i-1)} \!\geq\! q^\star$ by Lemma \ref{lemma:movement}. This implies $q^{(i)} \!\geq\! q^\star \!-\! \delta^{(i-1)}$ by the step-limiting constraint at $i \!-\! 1$. Since $q' \!=\! q^{(i)}$ and $\delta^{(i-1)} \!\leq\! \gamma^{-1}\delta^{(i)}$, this proves the lemma.\hfill$\square$

\begin{lemma}[Arbitrarily small $\delta$] \label{lemma:small}
For any tolerance $\varepsilon \!>\! 0$, there exists $K$ indicating the number of finite iterations, such that $\delta^{(K)} \!\leq\! \varepsilon$.
\end{lemma}
\textbf{Proof:}
Let $m^{(i)}$ denote the cumulative number of times the system has been in the oscillatory mode at iteration $i$, with $m^{(1)} \!=\! 0$ and $m^{(i+1)} \!=\! m^{(i)} + o^{(i)}$. Thus $\delta^{(i)} \!=\! \delta^{(0)}\gamma^{m^{(i)}}$. Following Lemma \ref{lemma:oscillatory}, for any $m \!>\! 0$, if $m^{(i)} \!=\! m$, there exists $l \!>\! 0$ such that $m^{(i+l)} \!=\! m+1$. Thus, we can make $m^{(i)}$ arbitrarily large with sufficient iterations, and therefore, $\delta^{(i)} \!=\! \delta^{(0)}\gamma^{m^{(i)}}$ can be made arbitrarily small. \hfill $\square$

\begin{figure*}[t]
    \centering
    \includegraphics[width=0.99\textwidth]{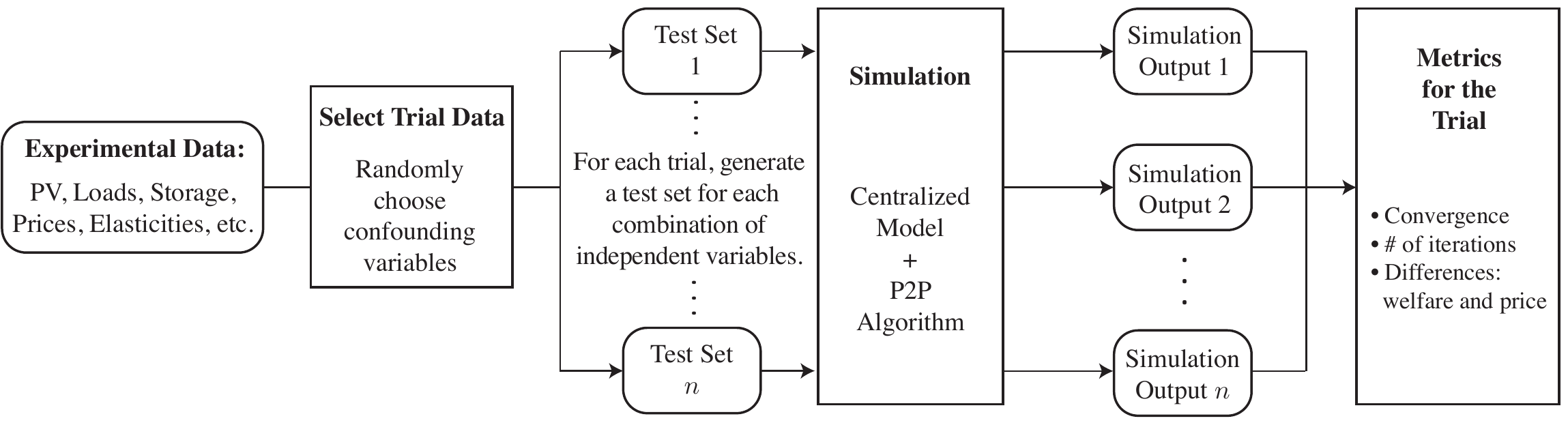}
    \caption{Trial procedure flowchart.}
    \label{fig:flowchart}
\end{figure*}

\begin{lemma}[Termination] \label{lemma:termination}
 If the algorithm terminates at iteration $i$ due to the stopping criterion, then $|q^{\prime,(i)} \!-\! q^\star| < \varepsilon$.
\end{lemma}
\textbf{Proof:}
We will prove the case when $q^\star \!<\! p_v^s$: First, consider the case when $q^{(i-1)} \!<\! q^\star$ (and hence $q^{\prime, (i-1)} \!=\! q^{(i-1)}$), then it follows from Lemma \ref{lemma:movement} that $q^{(i)} \!\ge\! q^{(i-1)}$. There are two cases, if (i) $q^{(i)} \!\ge\! q^\star$, then $q^\star \!<\! q^{\prime, (i)} \!\le\! q^{(i)}$ and it follows directly that $|q^{\prime, (i)} \!-\! q^\star| \!\le\! |q^{\prime, (i)} \!-\! q^{\prime, (i-1)}| \!\le\! \gamma \varepsilon \!<\! \varepsilon$, since the algorithm terminated at iteration $i$. Alternatively, if (ii) $q^{(i)} \!\le\! q^\star$, we have $q^{\prime, (i)} \!=\! q^{(i)}$ From \eqref{eq:q_sol} and the intersection of $g_k$ and $g_v$, it must be true that $q^{(i)} \!= q^{(i-1)} \!+\! \delta^{(i-1)}$, and thus $|q^{(i)} \!-\! q^{(i-1)}| = \delta^{(i-1)} \!\le\! \gamma \varepsilon$. Now, we have two cases: (a) if the system was oscillating, we have that $q^{\prime,(i-2)} \!>\! q^\star$ and so $q^{\prime,(i-2)} \!>\! q^{\prime,(i-1)}$ from Lemma 1. It is also true that $\delta^{(i-1)} \!=\! \gamma \delta^{(i-2)}$ and that $|q^{\prime,(i-1)} \!-\! q^{\prime,(i-2)}| \!\le\! \delta^{(i-2)}$, therefore, $|q^{\prime,(i)} \!-\! q^\star| \!\le\! |q^{\prime,(i)} \!-\! q^{\prime,(i-2)}| \!\le\! (1-\gamma)\delta^{(i-2)} \!=\! (1 \!-\! \gamma)\gamma^{-1}\delta^{(i-1)} \!\le\! (1 \!-\! \gamma)\gamma^{-1} \gamma \varepsilon \!=\! (1 \!-\! \gamma) \varepsilon \!<\! \varepsilon$. Or (b), if the system was not oscillating, this implies that $\delta^{(i-1)} \!=\! \delta^{(i-2)}$, and from the same argument as before satisfying that $|q^{(i)} \!-\! q^{(i-1)}| \!\le\! \gamma \varepsilon$. This is a contradiction, since that would also imply that $|q^{\prime, (i-1)} \!-\! q^{\prime, (i-2)}| \!\le\! \delta^{(i-2)} \!=\! \delta^{(i-1)} \!\le\! \gamma\varepsilon$, hence terminating before $i$. Second, for the case when $q^{(i-1)} \!>\! q^\star$, the proof is equivalent to the first case, but considering the special instance that if $q^{(i-1)} \!>\! p_v^s$, then $q^{\prime,(i-1)} \!=\! p_v^s$, but is still larger than $q^\star$ so the same idea holds by invoking Lemma \ref{lemma:movement}.

Finally, we address the case when $q^\star \!=\! p_v^s$. As described in Lemma \ref{lemma:oscillatory}, the algorithm will get stuck at $q' \!=\! p_v^s$ for more than two iterations, terminating the algorithm. Since $p_v^s \!=\! q^\star$, then $|q^{\prime,(i)} \!-\! q^{\prime,(i-1)}| \!=\! |q^{\prime,(i)} \!-\! q^\star| \!=\! 0 \!<\! \varepsilon$. \hfill $\square$

\begin{theorem}[Optimality of Algorithm 1]
\label{theorem1}
For 2 agents with strictly concave utility functions, $T=1$, and with sufficiently large max iterations $M$, Algorithm 1 returns a quantity within $\varepsilon$ of the centralized optimum $q^\star$.
\end{theorem}
\textbf{Proof:}
By Lemma \ref{lemma:small}, if we set $M \!\geq\! K$, then the algorithm will terminate due to the stopping criterion in at most $K$ iterations, and by Lemma \ref{lemma:termination}, the quantity is within $\varepsilon$ of $q^\star$. \hfill $\square$ 

In Fig.~\ref{fig:noncon} we depict a case that converges to the centralized solution via the step-limiting constraint. This case otherwise diverges based on the classic result \cite{kaldor1934classificatory} without the step-limiting constraint. This theoretical analysis provides the foundation for extending the algorithm to multiple agents $|\mathcal{C}| \!>\! 2$ with finite time horizon $T \!>\! 1$. We explore the behavior of the algorithm numerically for such cases in the next section.

\section{Computational Experiments and Simulations}

To provide additional insight into the algorithm performance, we perform two simulation-based computational experiments following the methodology and nomenclature in \cite{lara2020computational}. The simulation flowchart for both experiments are summarized in Fig. \ref{fig:flowchart}. The first, examines how the two algorithm parameters $\gamma$, $\delta^{(0)}$ affect the rate of convergence. The second, tests convergence for the unproven cases for $|\mathcal{C}| \!>\! 2$ and $T \!>\! 1$, and studies the effect of battery energy and power capacity on convergence and explores welfare differences between the centralized and P2P approaches.\footnote{The experiment parameters, data files, and MATLAB code to reproduce the experiments can all be found at https://github.com/Energy-MAC/TSG-P2P-Pricing.} In all experiments we use hourly load and PV profiles from Pecan Street \cite{pecan2020dataport}, and constant price elasticity utility functions fit to the baseline load with elasticities random on $[-1.5, -0.5]$ as in section \ref{subsec:centralized_simulation}.

\subsection{Effect of parameters $\gamma$ and $\delta^{(0)}$ on convergence}
\label{subsec:convergence}

In this experiment, we study the convergence rate for the $2$-agent, single period case. We systematically vary $\gamma \!\in\! \{0.05, 0.1, \dots, 0.95\}$, $\delta^{(0)} \!\in\! \{0.1, 0.2, \dots, 2\}$ kWh as independent variables, generating $380$ unique pairs of $(\gamma,\delta^{(0)})$. For each pair, we execute $100$ trials with randomly generated confounding variables (the two load profiles, hour of the year, price elasticities, and solar power between zero and twice the load) and compute the iterations to convergence. We use a stopping tolerance $\varepsilon = 10^{-3}$ for all trials.

\begin{figure}[t]
    \centering
\tikzset{every picture/.style={scale=1.0}}%
      \input{parameter_experiment}
  \caption{Effect of $\gamma$ on the number of iterations to convergence. The solid lines show the mean over all trials where $\delta^{(0)}$ lies in the interval specified in the legend. The dashed lines show the maximum.}
        \label{fig:parameter}
\end{figure}
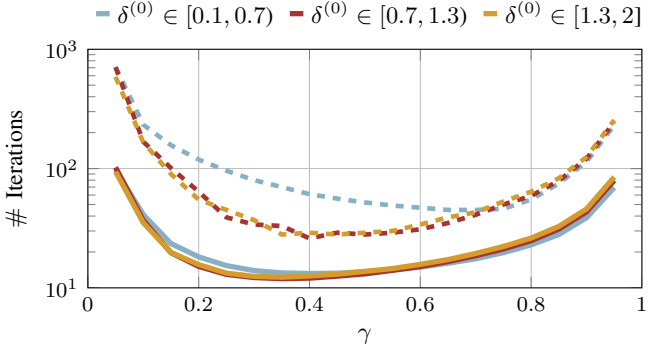

The results in Fig.~\ref{fig:parameter} show that $\gamma$ has a strong effect on the convergence rate and exhibits a minimum for $\gamma \!\in\! [0.3, 0.5]$ that is consistent across the different ranges of $\delta^{(0)}$, and that the algorithm converges in on the order of $10$-$20$ iterations on average for $\gamma$ in the middle range. We found that $\delta^{(0)}$ was not very significant in influencing the number of iterations except for causing an increase at especially small values, suggesting the  parameter ought to be set to a relatively large value. A possible intuition behind the effect of $\gamma$ is that especially small values shrink the box too quickly away from the equilibrium, while large values do not shrink rapidly enough.

\subsection{Performance for unproven cases}
To study the performance in the general (multi-agent) case, we vary the total battery capacity $\bar{S}_{\text{tot}}$\,$\in$\,$\{15,25,40,80,300\}$\,kWh and the maximum rate of charge/discharge of the battery $\bar{P}_{i}^b$\ $\in$\,$\{1,2,4,8\}$\,kW as independent variables, yielding $20$ distinct pairs.

Similar to section \ref{subsec:convergence}, for each pair $(\bar{S}_{\text{tot}}, \bar{P}_i^b)$ we execute $60$ trials ($60\times 20 \!=\! 1200$ simulations), randomly selecting PV and load profiles, price elasticities, an hour of the year, $T \!\in\! \{1,12,24\}$ hours, and number of agents $N \!\in\! [2,10]$. A battery capacity fraction is assigned uniformly to each agent (and then normalized) from the total battery capacity. The PV profiles are scaled so the total PV energy equals the total baseline load energy, and $(\gamma,\delta^{(0)}) \!=\! (0.5,0.5)$.

\subsubsection{Convergence performance}
All of the $1200$ treatments converge to a solution. The average iterations required to convergence is $112.5$, with a standard deviation of $257.3$ and a median of $61$. We observe that larger time horizons with more agents require more iterations for the algorithm to converge.

\subsubsection{Effect of battery parameters}
\begin{figure}[htbp]
  \centering
\tikzset{every picture/.style={scale=1.0}}%
\input{test_boxplot.tex}
 \caption{Number of iterations to converge with varying battery capacity.}
\label{fig:boxplot_size}
\end{figure}
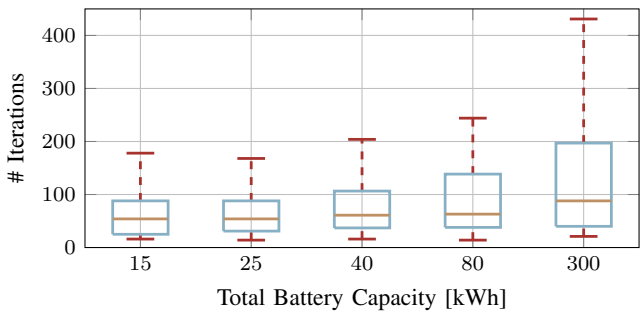

The effect of battery capacity on convergence is illustrated via boxplots in Fig.~\ref{fig:boxplot_size}, depicting the distribution of the number of iterations for convergence against battery capacity (with outliers omitted). 
In general, a higher battery capacity requires more iterations to converge. The intuition being that with higher battery availability, the flexibility for each agent to adapt to successive trades increases, thus requiring more iterations. This highlights the importance of storage in a P2P setting and the effect on the implementation of energy trading algorithms.
In contrast, the maximum charge/discharge rate of the battery does not significantly affect the number of iterations. This is expected, because given the demand profiles, a maximum rate of $1$ kW is usually enough to achieve a trade.

\subsubsection{Welfare comparison}
In order to compare the total welfare of all agents for the centralized and the iterative P2P algorithm, we classify the trials by grouping the time horizon. The statistics of welfare difference percentages $\Delta W_p$ and absolute welfare differences $\Delta W$ are presented in Table~\ref{tab:welf_diff}. We note that most of the entries for $\Delta W_p$ are lower than $0.1$\%, i.e., in the range of numerical tolerance used for MATLAB based optimizers.  These results indicate that in most cases the centralized welfare is close to that of the proposed algorithm. However, there exist cases when $T>1$, for which although the algorithm converges, the welfare is significantly different from the centralized solution.

\begin{table}[t]
\renewcommand{\arraystretch}{1.1}
    \caption{Welfare difference statistics for different time horizons.}
    \label{tab:welf_diff}
    \centering
    \begin{tabular}{l|l|rrr}
         & $T$ & $1$ & $12$ & $24$\\ \hline
         & \#Simulations &  $300$ & $420$ & $476$ \\ \hline
        & Mean [\%] & $0.023$ & $0.001$ & $0.072$ \\
       $\Delta W_p$ & Std [\%] & $0.079$ & $0.004$ & $0.706$ \\
        & Max [\%] & $0.558$ & $0.034$ & $7.758$ \\ \hline
         & Mean [\$] & $0.004$ & $0.002$ & $0.319$ \\
       $\Delta W$ & Std [\$] & $0.015$ & $0.007$ & $3.338$ \\
        & Max [\$] & $0.086$ & $0.057$ & $36.717$ \\ 
    \end{tabular}
\end{table}

\begin{table}[t]
\renewcommand{\arraystretch}{1.1}
    \caption{Welfare difference statistics for the special instance considered in Section IV-B.4.}
    \label{tab:welf_diff_case}
    \centering
    \begin{tabular}{l|r|rr|rr}
         & $W_\text{no}$ [\$] & $W_\text{centr}$ [\$]  & $W_\text{P2P}$ [\$]  & $\Delta W$ [\$] & $\Delta W_p$ [\%]\\ \hline
        ag-$1$  & $11.923$ & $19.804$ & $14.292$ & $5.512$ & $27.833$ \\
        ag-$2$  & $6.617$ & $16.785$ & $9.079$ & $7.706$ & $45.910$ \\
        ag-$3$  & $2.784$ & $2.933$ & $3.516$ & $-0.583$ & $-19.877$ \\
        ag-$4$  & $202.124$ & $202.906$ & $202.920$ & $-0.014$ & $-0.007$\\
        ag-$5$  & $164.184$ & $229.159$ & $203.345$ & $25.814$ & $11.265$\\
        $\pi$-ag & $1.633$ & $1.711$ & $3.429$ & $-1.718$ & $-100.409$ \\ \hline
        \textbf{Total} & $389.265$ & $473.298$ & $436.581$ & $36.717$ & $7.758$ \\
    \end{tabular}
\end{table}

\subsubsection{Special instance}
\begin{figure}[t]
    \centering
\tikzset{every picture/.style={scale=1.0}}%
      \input{bidding_plots.tex}
  \caption{Centralized and P2P algorithm price profiles for the special instance of considerable difference in welfare.}
        \label{figure:price_bidding}
\end{figure}
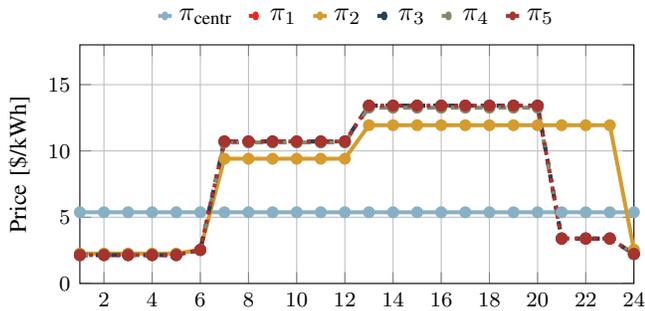

In this section we explore one instance where there is a considerable mismatch ($\Delta W \!=\! \$36.71$) between the welfare values obtained from the two approaches. This occurs for $T \!=\! 24$, $N \!=\! 6$, and low total battery capacity of $\bar{S}_{\text{tot}} \!=\! 15$ kWh. The key difference is that the prices for the agents in the algorithm are significantly different than those obtained in the centralized solution, as observed in Fig.~\ref{figure:price_bidding}. 

This simulation converges in $59$ iterations, when agent-$1$ (ag-$1$) exits the algorithm. However, at iteration $32$, agent-$2$ exits based on its stopping criteria, while the remaining agents continue trading, before exiting at iterations $59$, $58$, $58$, and $55$ respectively with similar price profiles, as indicated in Fig.~\ref{figure:price_bidding}. The consumption profiles and hence the individual welfare of each agent are thus significantly different from  the centralized solution. Table \ref{tab:welf_diff_case} summarizes the total welfare (consumption + trading) of each agent using the centralized and P2P algorithm. The welfare for the no-trading case $W_\text{no}$, is also presented for comparison.
A closer inspection reveals that while that agent-$3$ and the $\pi$-agent are better off in the P2P case, agents $1$, $2$, and $5$ are well placed in the centralized case. Furthermore, for this particular simulation, exiting earlier is not optimal for agent-$2$, although the price is lower than the other $q$-agents.

\subsubsection*{Summary} The simulation results in this section highlight the main contributions of this work:
\begin{enumerate}[label={(\arabic*)}, leftmargin=0.25cm, itemindent=0.3cm]
\setlength{\itemsep}{2pt}
    \item The P2P algorithm achieves similar welfare results as the centralized approach in most of the cases (Table~\ref{tab:welf_diff}), with the caveat that an early exit by some agents may introduce sub-optimality, in which case different agents end up as winners and losers relative to the social optimum (Table~\ref{tab:welf_diff_case}), but all agents are better off than no trading.
    \item More flexibility for the agents via larger storage or longer time horizons increases the number of iterations (Fig.~\ref{fig:boxplot_size}).
    \item The expected number of iterations is minimized by setting the shrinking rate of the step-size $\gamma$ around $0.4$ (Fig.~\ref{fig:parameter}).
    \item Real-time prices in a zero marginal-cost system arise from the marginal utility of consumption under scarcity.
\end{enumerate}

\section{Conclusions}
In this paper, we address the question of optimal pricing and mechanisms for achieving optimal dispatch in microgrids with scarce, zero marginal-cost energy resources. We contribute a novel analysis of a centralized economic dispatch with welfare maximization that uses a Lagrangian dual decomposition to state the equilibirum optimal price conditions and show a previously unstated result that although optimal prices can induce unique and optimal consumption profiles and generator output, they do not yield unique or power-balanced battery dispatch decisions except in particular circumstances. Next, we propose a P2P algorithm where agents keep utility functions private and iteratively interact by exchanging price and quantity offers to arrive at mutually agreeable and weakly Pareto-optimal trades. We theoretically prove this outcome converges to the social optimum within a specified tolerance for the $2$-agent case, and show via numerical experiments that the P2P algorithm converges in the multi-agent case, but we do not derive specific bounds. Although we find that the P2P algorithm obtains total welfare on average within $0.1$\% of the centralized solution for a wide range of parameters, significant differences in welfare and allocation can arise for longer time horizons and larger numbers of agents. We also find from simulations that the number of iterations for the P2P algorithm to converge increases with the total storage capacity, and that the P2P algorithm shrinking parameter $\gamma$ impacts the number of iterations, and should be set to the neighborhood of $0.3$-$0.5$ to minimize iterations in contexts similar to our simulations, while the initial maximum step-size $\delta^{(0)}$ is not significant.

The proposed P2P algorithm was designed to resemble an \textit{informal} decentralized trading process where prices arise from the value of electricity consumption under scarcity. We envision it is feasible to implement such an interaction in practice via a software platform that defines the rules and aids in the iteration, or even with informal negotiation between neighbors in a community. However, there are several limitations that need to be addressed for this approach to be useful in practice. First, we do not study the impact of strategic gaming between agents, which could be significant in small markets, nor the equity of outcomes. Conducting this analysis likely requires removing the assumption that $\pi$-agents offer prices equal to their dual variables and considering their profit maximizing strategy, given expectations of $q$-agents' demand curves. Second, our analysis only considers strictly concave utility functions. This is a common assumption, but may not capture the discrete nature of decisions to use particular loads at small time-scales. We expect it will be difficult for researchers to derive useful theoretical insights with non-concave utility functions, but the construction of realistic utility functions and consumption decision models for use in simulation would be of tremendous value to this and related work. Third, we do not include network constraints or validate power flow. While it is relatively straightforward to validate whether a particular negotiated dispatch is feasible given a network model, the impact of binding constraints on pricing and negotiation is non-trivial and warrants further study. Fourth, a system for matching agents into smaller negotiation pools based on expected outcomes may be necessary to handle large numbers of agents, e.g., hundreds. Here, a challenge is to design suitable exit strategies for satisfied agents without compromising the inviolability of agreements, and should also account for network constraints in creating market power (see \cite{kim2019p2p}). This introduces significant complexity, where methods to certify optimality or bound the outcome are important theoretical directions for future work. Lastly, we suggest the inclusion of uncertainty via scenarios in a stochastic programming framework to deal with uncertainty in solar forecasts and load estimation. The inclusion of power flow feasibility and network validation in the P2P algorithm, and extensions to a broader class of DERs are next steps in this research.

\appendices
\section{Non-ideal battery modeling}
\label{sec:detailed_battery_model}
In this appendix, we extend the analysis in Section II to a more realistic model of the battery and show that the results hold when battery inefficiencies and asymmetric charge and discharge power constraints are accounted for. We emphasize that we assume a linear model and the results may not hold for nonlinear models. However, agents may use decision models of varying complexity in practice, so understanding the implications of simplified models remains highly relevant.

The extended model replaces the net battery discharge $p^b_{i,t}$ with its positive discharge and charge components $p^{b,+}_{i,t}$ and $p^{b,-}_{i,t}$. We allow different discharge and charge power constraints $P_{i,t}^{b,+} \!>\! 0$ and $P_{i,t}^{b,-} \!>\! 0$ and assume that power is converted to and from the stored energy with charge efficiency $\sigma_i^- \in (0,1]$ and discharge efficiency $\sigma_i^+ \in (0,1]$, and that a battery self-discharges at a rate $(1\!-\!\theta_i)$ proportional to the state-of-charge, with $\theta_i \!\in\! [0,1)$. Thus, with the extended model, the problem formulation can be stated as:
\begin{subequations}
\label{eq:centralized_2}
\begin{align}
    &\min_{\boldsymbol{p}, \boldsymbol{d}, \boldsymbol{s}} ~ - \sum_{t\in \mathcal{T}} \sum_{n\in \mathcal{C}} U_{n,t} (d_{n,t}) \\
    &\text{ s.t. } ~~ \pi_t: \sum_{n\in \mathcal{C}} d_{n,t} = \sum_{i \in \mathcal{B}} (p_{i,t}^{b,+}-p_{i,t}^{b,-}) + \sum_{g\in \mathcal{G}} p_{g,t}^s, \, \forall t\in \mathcal{T} \label{eq:power_balance_2} \\
    &\hspace{0.9cm} \lambda_{g,t}^{s}: 0 \le p_{g,t}^{s} \le 
\bar{P}_{g,t}^s, \, \forall g \in \mathcal{G}, \forall t \in \mathcal{T} \label{eq:pv_lims_2}\\
&\hspace{0.9cm} {\lambda_{n,t}^{d,-}:-d_{n,t} \le 0, \, \forall n \in \mathcal{C}, \forall t \in \mathcal{T} \label{eq:demand_pos_2}}\\
&\hspace{0.9cm} \lambda_{i,t}^{b,+}:0 \le p_{i,t}^{b,+} \le \bar{P}_{i,t}^{b,+}, \, \forall i \in \mathcal{B}, \forall t \in \mathcal{T} \label{eq:batt_rate_pos_2}\\
&\hspace{0.9cm} \lambda_{i,t}^{b,-}:0 \le p_{i,t}^{b,-} \le \bar{P}_{i,t}^{b,-}, \, \forall i \in \mathcal{B}, \forall t \in \mathcal{T} \label{eq:batt_rate_neg_2}\\
&\hspace{1.9cm} p_{i,t}^{b,+} p_{i,t}^{b,-} = 0, \, \forall i \in \mathcal{B}, \forall t \in \mathcal{T} \label{eq:eq_batt_rate_comp_2}\\
&\hspace{0.9cm} \lambda_{i,t}^{c}:  0 \le s_{i,t}\le 
\bar{S}_{i,t}, \, \forall {i} \in \mathcal{B}, \forall t \in \mathcal{T} \label{eq:batt_lims_2} \\
&\hspace{1.9cm} s_{i,t} = \theta\, s_{i,t-1} + \sigma_i^-\, p_{i,t}^{b,-} -\notag\\ &\hspace{2.7cm} (\sigma_i^-)^{-1}p_{i,t}^{b,+} \Delta T, \, \forall {i}\in \mathcal{B}, \forall t \in \mathcal{T}. \label{eq:SOC_2}
\end{align}
\end{subequations}
In the following, we show a technical detail needed to extend \textbf{Remark 1}, state the modified forms of \textbf{Remark 2} and \textbf{Remark 4}, and show a modified example with discussion for \textbf{Remark 5}. \textbf{Remark 3} does not depend on the battery model and is unaffected.

For \textbf{Remark 1}, the analogous versions of \eqref{eq:W_nt}-\eqref{eq:dual_prob} follow mechanically in the same way as the main text in Section II. However, the complementarity constraint \eqref{eq:eq_batt_rate_comp_2} implies that $W_i(\pi)$ is now nonlinear, and thus strong duality does not necessarily hold, so an additional technical step is needed. We define a relaxed problem $W_i'(\pi)$ that drops the complementarity constraint to become a linear program. Under the assumption that $\pi_t \!\ge\! 0, \, \forall t$, which holds in particular for $\pi^\star$ by \textbf{Remark 2}, we show below that the problems have equal optimal objectives $W_i'(\pi) \!=\! W_i(\pi)$. Therefore, we can substitute $W_i'(\pi)$ for $W_i(\pi)$ into \eqref{eq:dual_prob}, and again rely on strong duality in $W_i'(\pi)$ to establish $\pi^\star$ as the optimal price.

To show that $W_i'(\pi)\!=\!W_i(\pi)$ if $\pi_t \!\ge\! 0, \, \forall t$, we take $(p'^{,b,+}_i,\,p'^{,b,\!-\!}_i)$ to be any any optimal solution to $W_i'(\pi)$, and compute a particular corresponding solution $(p^{\star,b,+}_i,\,p^{\star,b,\!-\!}_i)$. We verify below that $(p^{\star,b,+}_i,\,p^{\star,b,-}_i)$ is both feasible for the original problem and optimal for the relaxed problem, and therefore it is optimal for the original problem and each problem's optimal objectives are equal. The alternative solution is defined specifically to be the solution that has equal net charge/discharge as $(p'^{,b,+}_i,\,p'^{,b,\!-\!}_i)$ while satisfying the complementarity constraint: 
\begin{subequations}
\begin{align}
    p^{\star,b,+}_{i,t} &= \sigma_i^+\max((\sigma_i^+)^{-1}\,p'^{,b,+}_{i,t} \!-\! \sigma_i^-\, p'^{,b,-}_{i,t},\,0) \\
    p^{\star,b,-}_{i,t} &= (\sigma_i^-)^{-1}\max(\sigma_i^-\,p'^{,b,-}_{i,t} \!-\! (\sigma_i^+)^{-1}\,p'^{,b,+}_{i,t},\,0).
\end{align}
\end{subequations}
By construction, $ (\sigma_i^-)^{-1}\,p^{\star,b,+}_{i,t}-\sigma_i^-\, p^{\star,b,-}_{i,t} \!=\! (\sigma_i^-)^{-1}\,p'^{,b,+}_{i,t}-\sigma_i^-\, p'^{,b,-}_{i,t}$, and thus the stored energy trajectory is identical and \eqref{eq:batt_lims_2}-\eqref{eq:SOC_2} are satisfied. The power constraints \eqref{eq:batt_rate_pos_2}-\eqref{eq:batt_rate_neg_2} are satisfied because $0 \!\leq\! p^{\star,b,+}_{i,t} \!\leq\! p'^{,b,+}_{i,t}$ and  $0 \!\leq\! p^{\star,b,-}_{i,t} \!\leq\! p'^{,b,-}_{i,t}$. The complementarity constraint \eqref{eq:eq_batt_rate_comp_2} can be checked by inspection of the different cases of the max functions. Finally, $p^{\star,b,+}_{i,t}\!-\!p^{\star,b,-}_{i,t} \!\geq\! p'^{,b,+}_{i,t}\!-\!p'^{,b,-}_{i,t}, \, \forall t$, so if $\pi_t \!\geq\! 0$, then $-\sum_t \pi_t (p^{\star,b,+}_{i,t}-p^{\star,b,-}_{i,t}) \!\leq\! -\sum_t \pi_t (p'^{,b,+}_{i,t}-p'^{,b,-}_{i,t})$ and $(p^{\star,b,+},p^{\star,b,-})$ is optimal for the relaxed problem. \hfill $\square$

As an aside, note that if $\sigma_i^+\sigma_i^- \!<\! 1$ and $\pi_t \!>\! 0, \forall t$, then the relaxation is exact. This is because $(p^{\star,b,+},\,p^{\star,b,-}) \!\neq\! (p'^{,b,+}, p'^{,b,-})$ implies $\sum_t \pi_t (p^{\star,b,+}_{i,t}-p^{\star,b,-}_{i,t}) > \sum_t \pi_t (p'^{,b,+}_{i,t}-p'^{,b,-}_{i,t})$, which is a contradiction, therefore $(p^{\star,b,+}_i,\,p^{\star,b,-}_i) \!= \!(p'^{,b,+}_i,\, p'^{,b,-}_i)$.

The analogous form of \textbf{Remark 2} follows from the stationarity conditions of \eqref{eq:centralized_2}, and is
\begin{subequations}
\begin{align}
    \pi^\star_t& = {\partial U_{n,t}(d^\star_{n,t})}/{\partial d_{n,t}}+\lambda_{n, t}^{\star,d,-}, \, \forall n \in \mathcal{C} \label{eq:statio_d_2}\\
&=\lambda_{i,t}^{\star,b,+}\!-\! \Delta T\,(\sigma^+_i)^{-1}\,\sum_{\tau\geq t}\,\theta^{\tau-t}\!\,\lambda_{i,\tau}^{\star,c}, \, \forall i \in \mathcal{B} \label{eq:statio_batt_pos_2}\\
&=-\lambda_{i,t}^{\star,b,-}\!-\! \Delta T\,\sigma^-_i\,\sum_{\tau\geq t}\,\theta^{\tau-t}\!\,\lambda_{i,\tau}^{\star,c}, \, \forall i \in \mathcal{B} \label{eq:statio_batt_neg_2}\\
&=\lambda_{g,t}^{\star,s}, \, \forall g \in \mathcal{G}.
\end{align}
\end{subequations}

\textbf{Remark 4} again follows from \textbf{Remark 2}, and takes the form
\begin{subequations}
\begin{align}
\label{eq:rem4_2}
\theta\,\pi^\star_{t+1}-\pi^\star_t &= \theta\,\lambda_{i,t+1}^{\star,b,+}-\lambda_{i,t}^{\star,b,+}+\Delta T\,(\sigma^+_i)^{-1}\, \lambda_{i,t}^{\star,c} \\
&= \lambda_{i,t}^{\star,b,-}-\theta\,\lambda_{i,t+1}^{\star,b,-}+\Delta T\,\sigma ^-_i \,\lambda_{i,t}^{\star,c}.
\label{eq:rem4_2_charge}
\end{align}
\end{subequations}
These dynamics imply that the equilibrium optimal price, and equivalently the marginal value of consumption at optimum, will evolve depending on the battery inefficiencies and which battery constraints are active at the optimum.

For \textbf{Remark 5}, consider again the example from the main text with $T\!=\!5$, $\Delta T\!=\!1$, $\bar{P}^{b,+}_{i,t} \!\equiv\! 3$,  $\bar{S}_{i,t}\!=\!10$, $s_{0,i}\!=\!5$. Take $\sigma_i^+\!=\!0.95$, $\sigma_i^-\!=\!0.9$, $\theta_i\!=\!0.98$, $\bar{P}^{b,-}_{i,t} \equiv 2$, and a modified price $\pi^\star\!=\![1,1.0204,2,3,1.2680]$. Analogous to the example in the main text, solutions $p^b \!=\! [-0.9587, -0.9587, 3, 3, 0]$, $p^b \!=\! [0, -1.8983, 3, 3, 0]$, and $p^b \!=\! [-1.5863, -1.5863, 3, 3, 1]$ are all optimal solutions with net benefit $-13.063$.

In this example, $\pi_2$ was constructed by noting that if an optimal dispatch has the battery charging and unconstrained in power and stored energy at both $t\!=\!1$ and $t\!=\!2$, then by \eqref{eq:rem4_2}, $\pi_{t+1}=\theta^{-1}\pi_t$; i.e., the equilibrium price trajectory must have these dynamics when the storage is charging between successive periods and is not constrained. It then follows by the same principle as described in the main text that charge power can be feasibly shifted from one period to another without affecting the cost, thus the storage dispatch is not unique. Here $\pi_5$ was chosen as $\pi_5 = (\sigma^+\sigma^-)^{-1}\,\theta^{-4}$ to show equilibrium conditions where a higher future price exactly compensates for the lost energy from charging. Additional discussion about the impact of storage inefficiencies on optimal pricing can be found in~\cite{hug2016integration}.

\section{Constant price-elasticity utility functions}
\label{sec:utility_function}

Here we describe the procedure used for developing sample utility functions from data, assuming constant price-elasticities, which is one of two common simple assumptions in pricing theory (the other being a linear demand curve / quadratic utility function). We emphasize that this was chosen for example purposes only, and that all of the analysis only assumes that utility functions are strictly concave, and would apply to logarithmic or quadratic utility functions as well.

Let the marginal utility of consumption be denoted by $g(d) = \partial U(d)/\partial d$. As described in the main text, the equilibrium price is equal to the marginal utility of consumption, i.e., $\pi^\star = g(d^\star)$. Then the demand function of price is given by $h$ (as the inverse of $g$), where $d = h(\pi):=g^{-1}(\pi)$. The price-elasticity, which we denote as $r(\pi)$ is defined as the ratio of the percentage change in quantity to the percentage change in price, and in general depends on the price
\begin{align}
r(\pi) &= \frac{\text{d}h(\pi)}{\text{d}\pi}\frac{\pi}{h(\pi)}.
\end{align}
A constant price-elasticity implies $r(\pi) \equiv \hat{r}$. The general family of demand functions with this property has the form $h(\pi) = a\pi^{\hat{r}}$
for some constant $a$. This can be fit to an empirical price and consumption pair $(\pi_0,\,d_0)$ by setting $h(\pi_0)=d_0$ and obtaining $a=d_0\,\pi_0^{-\hat{r}}$. Inverting this to marginal utility and integrating to utility, one obtains
\begin{subequations}
\begin{align}
    g(d) &= \pi_0\left(\frac{d}{d_0}\right)^\frac{1}{\hat{r}}, \label{eq:cpe_g} \\
    U(d) &= \frac{\hat{r}\,\pi_0\, d^{\frac{1}{\hat{r}}+1}}{(\hat{r}+1)\,d_0^\frac{1}{\hat{r}}}+c.
\end{align}
\end{subequations}

However, general downward-sloping demand curves imply $r<0$, thus $\lim_{d\rightarrow 0^+}g(d)=\infty$ and $\lim_{d\rightarrow 0^+}U(d)=-\infty$ which can be problematic for optimization solvers and is also an unrealistic extreme in practice. Thus, we modify the function to have ``quasi-constant'' price elasticity, by shifting the marginal utility curve to the left by a small $\delta > 0$ and compensating the exponent for the shift so that $r(\pi_0)=\hat{r}$. We also choose a $c$ such that $U(0)=0$. The resulting marginal utility and utility functions are
\begin{subequations}
\begin{align}
    g(d) &= \pi_0\left(\frac{d+\delta}{d_0+\delta}\right)^\frac{1}{r'}, \\
    U(d) &= \frac{r'\,\pi_0\left((d+\delta)^{\frac{1}{r'}+1}-\delta^{\frac{1}{r'}+1}\right)}{(r'+1)(d_0+\delta)^\frac{1}{r'}},\\
    r' &= \hat{r}\left(1+\frac{\delta}{d_0}\right)^{-1}.
\end{align}
\end{subequations}

\bibliographystyle{IEEEtran}
\input{main.bbl}

\end{document}

%% file: centralized_plots.tex
%
\definecolor{mycolor5}{RGB}{136,176,197}
\definecolor{mycolor1}{RGB}{168,50,45}
\definecolor{mycolor4}{HTML}{253F5B}%
\definecolor{mycolor2}{HTML}{818A6F}
\definecolor{mycolor3}{HTML}{D59B2D}
\begin{tikzpicture}

\begin{axis}[%
width=2.9in,
height=0.8in,
at={(1.394in,5.744in)},
scale only axis,
title = {\small (a) Total Solar Output [kW]},
title style={at={(0.5,1.1)},anchor=north,yshift=-0.1},
xmin=1,
xmax=66,
ymin=0,
ymax=40,
ylabel style={font=\color{black}},
yticklabel style = {font=\footnotesize,xshift=0ex},
xticklabel style = {font=\footnotesize,yshift=0ex},
axis background/.style={fill=white},
xmajorgrids,
ymajorgrids,
legend style={legend cell align=left, align=left, draw=black}
]
\addplot [color=mycolor1,line width=1.5pt]
  table[row sep=crcr]{%
1	2.0240011299435\\
2	3.4667927672956\\
3	6.21084999999999\\
4	12.7726333333333\\
5	22.7006333333333\\
6	32.067\\
7	30.42645\\
8	22.8280351851852\\
9	25.3884764705882\\
10	11.7821166666667\\
11	5.32863333333333\\
12	1.68163333333334\\
13	0\\
23	0\\
24	0.35456666666667\\
25	2.68778898305085\\
26	8.03556666666667\\
27	17.4883833333333\\
28	24.5638333333333\\
29	23.4538318181818\\
30	32.2766833333333\\
31	27.0679666666667\\
32	30.6587833333333\\
33	23.9247666666667\\
34	16.5855547619048\\
35	8.77723245614035\\
36	2.04741666666666\\
37	0\\
48	0\\
49	1.31704999999999\\
50	4.2111\\
51	7.21071666666667\\
52	13.3404761904762\\
53	16.3772166666667\\
54	17.0947166666667\\
55	15.9243\\
56	11.5906372881356\\
57	8.4539880952381\\
58	5.99516818181819\\
59	2.13211666666666\\
60	0.093141156462579\\
61	0\\
71	0\\
72	0.694933333333339\\
};

\addplot [color=mycolor2,line width=1.5pt]
 table[row sep=crcr]{%
1	2.0240011299435\\
2	3.4667927672956\\
3	6.21084999999999\\
4	12.7726333333333\\
5	22.7006333333333\\
6	32.067\\
7	30.42645\\
8	22.8280351851852\\
9	25.3884764705882\\
10	11.7821166666667\\
11	5.32863333333333\\
12	1.68163333333334\\
13	0\\
23	0\\
24	0.35456666666667\\
25	2.68778898305085\\
26	8.03556666666667\\
27	17.4883833333333\\
28	24.5638333333333\\
29	23.4538318181818\\
30	32.2766833333333\\
31	27.0679666666667\\
32	30.6587833333333\\
33	23.9247666666667\\
34	16.5855547619048\\
35	8.77723245614035\\
36	2.04741666666666\\
37	0\\
48	0\\
49	1.31704999999999\\
50	4.2111\\
51	7.21071666666667\\
52	13.3404761904762\\
53	16.3772166666667\\
54	17.0947166666667\\
55	15.9243\\
56	11.5906372881356\\
57	8.4539880952381\\
58	5.99516818181819\\
59	2.13211666666666\\
60	0.093141156462579\\
61	0\\
71	0\\
72	0.694933333333339\\
};

\addplot [color=mycolor3,line width=1.5pt]
  table[row sep=crcr]{%
1	2.0240011299435\\
2	3.4667927672956\\
3	6.21084999999999\\
4	12.7726333333333\\
5	22.7006333333333\\
6	32.067\\
7	30.42645\\
8	22.8280351851852\\
9	25.3884764705882\\
10	11.7821166666667\\
11	5.32863333333333\\
12	1.68163333333334\\
13	0\\
23	0\\
24	0.35456666666667\\
25	2.68778898305085\\
26	8.03556666666667\\
27	17.4883833333333\\
28	24.5638333333333\\
29	23.4538318181818\\
30	32.2766833333333\\
31	27.0679666666667\\
32	30.6587833333333\\
33	23.9247666666667\\
34	16.5855547619048\\
35	8.77723245614035\\
36	2.04741666666666\\
37	0\\
48	0\\
49	1.31704999999999\\
50	4.2111\\
51	7.21071666666667\\
52	13.3404761904762\\
53	16.3772166666667\\
54	17.0947166666667\\
55	15.9243\\
56	11.5906372881356\\
57	8.4539880952381\\
58	5.99516818181819\\
59	2.13211666666666\\
60	0.093141156462579\\
61	0\\
71	0\\
72	0.694933333333339\\
};

\addplot [color=mycolor4,line width=1.5pt]
  table[row sep=crcr]{%
1	2.0240011299435\\
2	3.4667927672956\\
3	6.21084999999999\\
4	12.7726333333333\\
5	22.7006333333333\\
6	32.067\\
7	30.42645\\
8	22.8280351851852\\
9	25.3884764705882\\
10	11.7821166666667\\
11	5.32863333333333\\
12	1.68163333333334\\
13	0\\
23	0\\
24	0.35456666666667\\
25	2.68778898305085\\
26	8.03556666666667\\
27	17.4883833333333\\
28	24.5638333333333\\
29	23.4538318181818\\
30	32.2766833333333\\
31	27.0679666666667\\
32	30.6587833333333\\
33	23.9247666666667\\
34	16.5855547619048\\
35	8.77723245614035\\
36	2.04741666666666\\
37	0\\
48	0\\
49	1.31704999999999\\
50	4.2111\\
51	7.21071666666667\\
52	13.3404761904762\\
53	16.3772166666667\\
54	17.0947166666667\\
55	15.9243\\
56	11.5906372881356\\
57	8.4539880952381\\
58	5.99516818181819\\
59	2.13211666666666\\
60	0.093141156462579\\
61	0\\
71	0\\
72	0.694933333333339\\
};

\end{axis}

\begin{axis}[%
width=2.9in,
height=0.8in,
at={(1.394in,4.644in)},
scale only axis,
title = {\small (b) Total Power Consumption [kW]},
title style={at={(0.5,1.1)},anchor=north,yshift=-0.1},
xmin=1,
xmax=66,
ymin=0,
ymax=45,
ylabel style={font=\color{black}},
yticklabel style = {font=\footnotesize,xshift=0ex},
xticklabel style = {font=\footnotesize,yshift=0ex},
axis background/.style={fill=white},
xmajorgrids,
ymajorgrids,
legend style={legend cell align=left, align=left, draw=black, font=\small, draw=none, legend columns=-1, at={(1.02,2.85)}}
]
\addplot [color=mycolor1,line width=1.5pt]
  table[row sep=crcr]{%
1	1.71243268947129\\
2	2.9880967756857\\
3	6.99960259312313\\
4	12.7679049512147\\
5	22.6986060061713\\
6	25.4817874769851\\
7	27.0248042315461\\
8	22.8236023313379\\
9	25.3885847108417\\
10	11.7822515055246\\
11	5.33008520446204\\
12	4.68116590459799\\
13	5.99243711226522\\
14	0.225176575755583\\
15	0.199221441879544\\
16	0.115706139438643\\
17	0.0851195513260308\\
18	0.064172356741949\\
20	0.0567588501368022\\
21	0.0838161171469523\\
23	0.0520160564842769\\
24	0.354596550446814\\
25	2.686304505103\\
26	7.81683060402082\\
27	17.6862628771424\\
28	17.3373681891172\\
29	24.7949303606195\\
30	28.5228481222328\\
31	26.7348425688347\\
32	30.6542363409397\\
33	23.9261878275201\\
34	16.5843138746375\\
35	8.77714227515733\\
36	5.09146407500184\\
37	5.64319876666336\\
38	0.276266543049132\\
39	0.263868463718538\\
40	0.153854046485776\\
41	0.0965039039736553\\
42	0.0985529725452921\\
44	0.0548338693350701\\
45	0.0836556163952196\\
47	0.0404001193756187\\
48	0.106208211608106\\
49	1.31676127339365\\
50	2.42082803660851\\
51	6.58876687165387\\
52	14.4972914522886\\
53	17.6164366646441\\
54	9.08624663057482\\
55	14.2175172621325\\
56	11.3230199496448\\
57	8.45423486129057\\
58	5.99530370709182\\
59	3.15486921306038\\
60	2.85761865260326\\
61	5.49165551223527\\
62	0.176974187998795\\
63	0.137595224692916\\
64	0.1063289555496\\
65	0.114075329403704\\
66	0.0277685928790419\\
67	0.0586295957215697\\
68	0.0188925340509343\\
69	0.0171387077285061\\
70	0.036796591202787\\
71	0.0258611505738173\\
72	0.694883304154857\\
};
\addlegendentry{$\bar{S}_{\text{tot}}\text{ = $1$}~~$}

\addplot [color=mycolor2,line width=1.5pt]
  table[row sep=crcr]{%
1	1.71216740566493\\
2	2.98789216106682\\
3	6.99962685588321\\
4	10.4904024317957\\
5	14.8440048744299\\
6	13.8367330015257\\
7	14.7192426226889\\
8	16.9068166079617\\
9	25.3890150249179\\
10	17.4980437838407\\
11	13.902604801838\\
12	14.3653158497719\\
13	18.5228658728713\\
14	0.84026667550097\\
15	0.758559314538829\\
16	0.502263274097572\\
17	0.409621450359111\\
18	0.338370780676527\\
19	0.30439515558821\\
20	0.308998682152321\\
21	0.396854329772708\\
23	0.280830585371319\\
24	0.376704306717826\\
25	2.68519217430672\\
26	5.72901461971551\\
27	12.9089886691259\\
28	11.7982225445873\\
29	16.8906977318402\\
30	19.4692358774663\\
31	18.3334499120168\\
32	28.4188446968832\\
33	23.9261782128255\\
34	16.5868130792216\\
35	15.8538647942203\\
36	18.1579890744493\\
37	20.2796061560219\\
38	1.18160953968611\\
39	1.13545093484925\\
40	0.738496014868829\\
41	0.526295519289462\\
42	0.514478166689997\\
43	0.439402540954504\\
44	0.352713212099019\\
45	0.468235221215735\\
46	0.374509877432743\\
47	0.272818127834569\\
48	0.525037994290813\\
49	1.3172523791914\\
50	1.161082244821\\
51	3.23148377007519\\
52	7.05542434420296\\
53	8.55251720386146\\
54	4.10944174769378\\
55	6.46786932955312\\
56	5.17274766185126\\
57	8.55364280090905\\
58	12.930924246101\\
59	11.2044299040996\\
60	10.1449793381615\\
61	19.6364095058984\\
62	0.813472759865547\\
63	0.685814895446939\\
64	0.548990690267303\\
65	0.577260550518432\\
66	0.251625020124848\\
67	0.358165510831228\\
68	0.217767292619982\\
69	0.217755923802869\\
70	0.289494305632161\\
71	0.242219184069782\\
72	0.695497821357662\\
};
\addlegendentry{$\bar{S}_{\text{tot}}\text{ = $5$}~~$}

\addplot [color=mycolor3,line width=1.5pt]
  table[row sep=crcr]{%
1	0.517994517964368\\
2	0.93916712158601\\
3	2.27838043481468\\
4	1.92512159862599\\
5	2.72597183809062\\
6	2.53858196706173\\
7	2.72777294494102\\
8	3.15220769800426\\
9	27.9499405779377\\
10	33.2252325741915\\
11	26.3459490743132\\
12	27.3007279084505\\
13	35.3794192724339\\
14	1.66467101857984\\
15	1.50251409465001\\
16	1.02518742501807\\
17	0.854258728556246\\
18	0.718670181632291\\
19	0.654465378352128\\
20	0.662964713626891\\
21	0.83811432832708\\
23	0.604777344452216\\
24	0.769735295473552\\
25	0.635995562495779\\
26	0.893965125330141\\
27	2.11302742988681\\
28	1.92413969829181\\
29	2.79088335489116\\
30	3.24485823318868\\
31	3.06078564080191\\
32	4.73786815814114\\
33	29.6395939573465\\
34	28.1821422040617\\
35	31.9938595954456\\
36	36.7891198600909\\
37	41.2846267375531\\
38	2.47401304868214\\
39	2.38315578148655\\
40	1.57883391515244\\
41	1.14963177726064\\
42	1.1215360413813\\
43	0.967665260810804\\
44	0.799011535762389\\
45	1.03715938914777\\
46	0.844033317731117\\
47	0.634889364345284\\
48	1.14158101057215\\
49	1.51073968216139\\
50	0.463851768589592\\
51	1.38118352848235\\
52	3.0269983997005\\
53	3.66877176475396\\
54	1.75757027589786\\
55	2.79864365681874\\
56	2.23358756763419\\
57	13.0561090314007\\
58	19.9070726462851\\
59	17.2383034745309\\
60	15.5988155026925\\
61	30.3317161364227\\
62	1.30143498076737\\
63	1.10300767546431\\
64	0.885859296590283\\
65	0.931261705420638\\
66	0.431520508552921\\
67	0.59527573822804\\
68	0.379227488415339\\
69	0.379521624048323\\
70	0.491818725050507\\
71	0.416347995611247\\
72	0.696745828174897\\
};
\addlegendentry{$\bar{S}_{\text{tot}}\text{ = $15$}~~$}

\addplot [color=mycolor4,line width=1.5pt]
  table[row sep=crcr]{%
1	0.517943053641446\\
2	0.939084192843538\\
3	2.2781714416757\\
4	1.92496703870005\\
5	2.72582418109222\\
6	2.53841234464858\\
7	2.72762845935634\\
8	3.15207923560411\\
9	27.9496547495341\\
10	33.2247892947707\\
11	26.3462191473873\\
12	27.3010138594626\\
13	35.3802052966422\\
14	1.66471359999743\\
15	1.50257179577268\\
16	1.02522361752571\\
17	0.854283230786905\\
18	0.718712007551247\\
19	0.654497639862058\\
20	0.663003489885767\\
21	0.838160699604998\\
23	0.604819258057745\\
24	0.769792611342595\\
25	0.635687750552052\\
26	0.893561183755835\\
27	2.11206347770563\\
28	1.92333702913359\\
29	2.78971604687102\\
30	3.24354011864719\\
31	3.05961034068615\\
32	4.73617030505859\\
33	29.6357441952206\\
34	28.180164298161\\
35	31.9926621693781\\
36	36.7887434142534\\
37	41.2846959837276\\
38	2.47398920684448\\
39	2.38316426121845\\
40	1.57884773562733\\
41	1.1496410175013\\
42	1.12154710802774\\
43	0.967683099809037\\
44	0.799037800760203\\
45	1.03719262663587\\
46	0.844067282686069\\
47	0.634925034557597\\
48	1.14164101243557\\
49	1.51084782627883\\
50	0.46388862556843\\
51	1.38127923766443\\
52	3.02720059031964\\
53	3.66895716147984\\
54	1.75769698859951\\
55	2.79885542820602\\
56	2.2337952905124\\
57	13.057498179135\\
58	19.9092994580971\\
59	17.2406249472676\\
60	15.6008953517159\\
61	30.3362558633935\\
62	1.30163276463226\\
63	1.10317768688897\\
64	0.886006816946107\\
65	0.931427952192109\\
66	0.431595542516376\\
67	0.595387895796151\\
68	0.379309315575355\\
69	0.379605199736687\\
70	0.491920693871847\\
71	0.416438719198155\\
72	0.696747068906063\\
};
\addlegendentry{$\bar{S}_{\text{tot}}\text{ = $300$ [kWh]}$}

\end{axis}

\begin{axis}[%
width=2.9in,
height=0.8in,
at={(1.394in,3.544in)},
scale only axis,
title = \small{(c) Price profile [\$]},
title style={at={(0.5,1.1)},anchor=north,yshift=-0.1},
xmin=1,
xmax=66,
ymin=0,
ymax=0.6,
ylabel style={font=\color{black}},
yticklabel style = {font=\footnotesize,xshift=0ex},
xticklabel style = {font=\footnotesize,yshift=0ex},
axis background/.style={fill=white},
xmajorgrids,
ymajorgrids,
legend style={legend cell align=left, align=left, draw=black}
]
\addplot [color=mycolor1,line width=1.5pt]
  table[row sep=crcr]{%
1	0.162920686228105\\
2	0.162721903220898\\
3	0.162647711764521\\
4	0.120321380754433\\
5	0.110071622124551\\
6	0.102336935187068\\
7	0.102346156372576\\
8	0.1154704102632\\
9	0.260487899416006\\
10	0.376941253647288\\
11	0.470449332888293\\
12	0.500401769553378\\
14	0.50035742004728\\
17	0.500048466190819\\
20	0.499451193537794\\
21	0.499141898414322\\
22	0.498609792972132\\
23	0.497394472958561\\
24	0.325510034758949\\
25	0.138625930246178\\
26	0.103600509810676\\
27	0.103317073091148\\
28	0.100501860104615\\
30	0.10052111972135\\
31	0.100704024608959\\
32	0.113521286383488\\
33	0.257771461303761\\
34	0.292084872428347\\
35	0.394654730166721\\
36	0.514555169537445\\
38	0.514525524317335\\
41	0.514298282622221\\
44	0.51385829144327\\
46	0.513339501065118\\
47	0.512910251768105\\
48	0.512315378253177\\
49	0.249283988044567\\
50	0.12899818595379\\
53	0.128908631330361\\
54	0.125086738191982\\
55	0.125111428675453\\
56	0.1253436308196\\
57	0.280921837101175\\
58	0.37816267684957\\
59	0.458797115938751\\
61	0.458859093414603\\
64	0.458515020188401\\
65	0.458307795940371\\
67	0.457569419362088\\
68	0.456929115907371\\
69	0.456136229692987\\
70	0.455125245742721\\
71	0.452681695932313\\
72	0.206470739165297\\
};

\addplot [color=mycolor2,line width=1.5pt]
  table[row sep=crcr]{%
1	0.162930495744178\\
2	0.162726427070837\\
3	0.162647209788403\\
4	0.129953126735117\\
8	0.129931621986671\\
9	0.260486277471728\\
10	0.322695280534518\\
16	0.322658063665955\\
21	0.322457932851123\\
23	0.32215042178558\\
24	0.319305883584065\\
25	0.138647681341936\\
26	0.116962275859052\\
31	0.116931267102842\\
32	0.116958352593969\\
33	0.257771526037615\\
34	0.292066829604806\\
35	0.312549954125714\\
41	0.312508151476763\\
46	0.312341278362823\\
48	0.312142153157438\\
49	0.249249626591208\\
50	0.170743103581955\\
56	0.170722424269471\\
57	0.279635735259134\\
58	0.280272935428982\\
64	0.280228751119921\\
68	0.280049578568025\\
70	0.279846361371895\\
71	0.279584651159993\\
72	0.206406064485478\\
};

\addplot [color=mycolor3,line width=1.5pt]
  table[row sep=crcr]{%
1	0.250919641748652\\
5	0.250852324797052\\
22	0.250709723565265\\
24	0.250496301205729\\
25	0.236885892574463\\
29	0.236846099864223\\
67	0.236733779264711\\
70	0.236577116426119\\
71	0.236424167868947\\
72	0.206276067608627\\
};

\addplot [color=mycolor4,line width=1.5pt]
  table[row sep=crcr]{%
1	0.250927450066655\\
5	0.250858106853059\\
21	0.250741190596216\\
23	0.250641443996628\\
24	0.250489299452241\\
25	0.236925515758131\\
29	0.236883208791681\\
65	0.236770396208968\\
70	0.236560461860037\\
71	0.236406921276838\\
72	0.206275107702751\\
};

\end{axis}

\begin{axis}[%
width=2.9in,
height=0.8in,
at={(1.394in,2.444in)},
scale only axis,
title = {\small (d) Total Battery Output [kW]},
title style={at={(0.5,1.1)},anchor=north,yshift=-0.1},
xmin=1,
xmax=66,
ymin=-35,
ymax=45,
yticklabel style = {font=\footnotesize,xshift=0ex},
xticklabel style = {font=\footnotesize,yshift=0ex},
ylabel style={font=\color{black}},
axis background/.style={fill=white},
xmajorgrids,
ymajorgrids,
legend style={legend cell align=left, align=left, draw=black}
]
\addplot [color=mycolor1,line width=1.5pt]
  table[row sep=crcr]{%
1	-0.311568440472229\\
2	-0.478695991609897\\
3	0.788752593123135\\
4	-0.00472838211864257\\
5	-0.00202732716199705\\
6	-6.58521252301493\\
7	-3.40164576845395\\
8	-0.0044328538473053\\
11	0.00145187112870815\\
13	5.99243711226522\\
14	0.225176575755583\\
15	0.199221441879544\\
16	0.115706139438643\\
18	0.064172356741949\\
20	0.0567588501368022\\
21	0.0838161171469523\\
23	0.0520160564842769\\
24	2.98837801580021e-05\\
25	-0.00148447794785511\\
26	-0.218736062645846\\
27	0.197879543809037\\
28	-7.22646514421612\\
29	1.34109854243772\\
30	-3.75383521110048\\
31	-0.333124097832012\\
32	-0.00454699239358547\\
34	-0.00124088726730065\\
35	-9.0180983022492e-05\\
36	3.04404740833516\\
37	5.64319876666336\\
38	0.276266543049132\\
39	0.263868463718538\\
40	0.153854046485776\\
41	0.0965039039736553\\
42	0.0985529725452921\\
44	0.0548338693350701\\
45	0.0836556163952196\\
47	0.0404001193756187\\
48	0.106208211608106\\
49	-0.000288726606342493\\
50	-1.79027196339148\\
51	-0.621949795012782\\
52	1.15681526181237\\
53	1.23921999797741\\
54	-8.00847003609184\\
55	-1.7067827378675\\
56	-0.267617338490751\\
57	0.000246766052470093\\
58	0.000135525273648796\\
59	1.02275254639372\\
60	2.76447749614067\\
61	5.49165551223527\\
62	0.176974187998795\\
64	0.1063289555496\\
65	0.114075329403704\\
66	0.0277685928790419\\
67	0.0586295957215697\\
68	0.0188925340509343\\
69	0.0171387077285061\\
70	0.036796591202787\\
71	0.0258611505738173\\
72	-5.00291784817364e-05\\
};

\addplot [color=mycolor2,line width=1.5pt]
  table[row sep=crcr]{%
1	-0.311833724278586\\
2	-0.478900606228777\\
3	0.788776855883214\\
4	-2.2822309015376\\
5	-7.85662845890344\\
6	-18.2302669984743\\
7	-15.7072073773111\\
8	-5.92121857722347\\
9	0.000538554329622798\\
10	5.71592711717402\\
11	8.57397146850469\\
12	12.6836825164386\\
13	18.5228658728713\\
14	0.84026667550097\\
15	0.758559314538829\\
16	0.502263274097572\\
17	0.409621450359111\\
18	0.338370780676527\\
19	0.30439515558821\\
20	0.308998682152321\\
21	0.396854329772708\\
23	0.280830585371319\\
24	0.0221376400511701\\
25	-0.00259680874412993\\
26	-2.30655204695115\\
27	-4.57939466420747\\
28	-12.765610788746\\
29	-6.56313408634163\\
30	-12.807447455867\\
31	-8.73451675464986\\
32	-2.23993863645011\\
33	0.00141154615886308\\
34	0.00125831731682524\\
35	7.07663233807999\\
36	16.1105724077827\\
37	20.2796061560219\\
38	1.18160953968611\\
39	1.13545093484925\\
40	0.738496014868829\\
41	0.526295519289462\\
42	0.514478166689997\\
43	0.439402540954504\\
44	0.352713212099019\\
45	0.468235221215735\\
47	0.272818127834569\\
48	0.525037994290813\\
49	0.000202379191406976\\
50	-3.050017755179\\
51	-3.97923289659148\\
52	-6.28505184627323\\
53	-7.82469946280521\\
54	-12.9852749189729\\
55	-9.45643067044688\\
56	-6.41788962628434\\
57	0.0996547056709431\\
58	6.93575606428284\\
59	9.07231323743289\\
60	10.0518381816989\\
61	19.6364095058984\\
62	0.813472759865547\\
64	0.548990690267303\\
65	0.577260550518432\\
66	0.251625020124848\\
67	0.358165510831228\\
68	0.217767292619982\\
69	0.217755923802869\\
70	0.289494305632161\\
71	0.242219184069782\\
72	0.00056448802432385\\
};

\addplot [color=mycolor3,line width=1.5pt]
  table[row sep=crcr]{%
1	-1.50600661197913\\
2	-2.52762564570959\\
3	-3.93246956518531\\
4	-10.8475117347073\\
5	-19.9746614952427\\
6	-29.5284180329383\\
7	-27.698677055059\\
8	-19.6758274871809\\
9	2.56146410734947\\
10	21.4431159075248\\
11	21.0173157409799\\
12	25.6190945751172\\
13	35.3794192724339\\
14	1.66467101857984\\
15	1.50251409465001\\
16	1.02518742501807\\
17	0.854258728556246\\
18	0.718670181632291\\
19	0.654465378352128\\
20	0.662964713626891\\
21	0.83811432832708\\
23	0.604777344452216\\
24	0.415168628806882\\
25	-2.05179342055507\\
26	-7.14160154133653\\
27	-15.3753559034465\\
28	-22.6396936350415\\
29	-20.6629484632906\\
30	-29.0318251001447\\
31	-24.0071810258647\\
32	-25.9209151751922\\
33	5.71482729067985\\
34	11.596587442157\\
36	34.7417031934243\\
37	41.2846267375531\\
38	2.47401304868214\\
39	2.38315578148655\\
40	1.57883391515244\\
41	1.14963177726064\\
42	1.1215360413813\\
43	0.967665260810804\\
44	0.799011535762389\\
45	1.03715938914777\\
46	0.844033317731117\\
47	0.634889364345284\\
48	1.14158101057215\\
49	0.193689682161391\\
50	-3.74724823141041\\
51	-5.8295331381843\\
52	-10.3134777907757\\
53	-12.7084449019127\\
54	-15.3371463907688\\
55	-13.1256563431813\\
56	-9.35704972050139\\
57	4.60212093616256\\
58	13.9119044644669\\
59	15.1061868078643\\
60	15.5056743462299\\
61	30.3317161364227\\
62	1.30143498076737\\
63	1.10300767546431\\
64	0.885859296590283\\
65	0.931261705420638\\
66	0.431520508552921\\
67	0.59527573822804\\
68	0.379227488415339\\
69	0.379521624048323\\
70	0.491818725050507\\
71	0.416347995611247\\
72	0.0018124948415732\\
};

\addplot [color=mycolor4,line width=1.5pt]
  table[row sep=crcr]{%
1	-1.50605807630207\\
2	-2.52770857445206\\
3	-3.9326785583243\\
4	-10.8476662946333\\
5	-19.9748091522411\\
6	-29.5285876553514\\
7	-27.6988215406436\\
8	-19.6759559495811\\
9	2.56117827894587\\
10	21.442672628104\\
11	21.017585814054\\
12	25.6193805261292\\
13	35.3802052966422\\
14	1.66471359999743\\
15	1.50257179577268\\
16	1.02522361752571\\
17	0.854283230786905\\
18	0.718712007551247\\
19	0.654497639862058\\
20	0.663003489885767\\
21	0.838160699604998\\
23	0.604819258057745\\
24	0.415225944675939\\
25	-2.0521012324988\\
26	-7.14200548291083\\
27	-15.3763198556277\\
28	-22.6404963041997\\
29	-20.6641157713108\\
30	-29.0331432146861\\
31	-24.0083563259805\\
32	-25.9226130282747\\
33	5.71097752855391\\
34	11.5946095362562\\
36	34.7413267475868\\
37	41.2846959837276\\
38	2.47398920684448\\
39	2.38316426121845\\
40	1.57884773562733\\
41	1.1496410175013\\
42	1.12154710802774\\
43	0.967683099809037\\
44	0.799037800760203\\
45	1.03719262663587\\
46	0.844067282686069\\
47	0.634925034557597\\
48	1.14164101243557\\
49	0.193797826278825\\
50	-3.74721137443157\\
51	-5.82943742900224\\
52	-10.3132756001566\\
53	-12.7082595051868\\
54	-15.3370196780672\\
55	-13.125444571794\\
56	-9.35684199762319\\
57	4.60351008389689\\
58	13.9141312762789\\
59	15.108508280601\\
60	15.5077541952533\\
61	30.3362558633935\\
62	1.30163276463226\\
63	1.10317768688897\\
64	0.886006816946107\\
65	0.931427952192109\\
66	0.431595542516376\\
67	0.595387895796151\\
68	0.379309315575355\\
69	0.379605199736687\\
70	0.491920693871847\\
71	0.416438719198155\\
72	0.00181373557273901\\
};
\end{axis}

\begin{axis}[%
width=2.9in,
height=0.8in,
at={(1.394in,1.344in)},
scale only axis,
title = {\small (e) Mean State-of-Charge [\%]},
title style={at={(0.5,1.1)},anchor=north,yshift=-0.1},
xmin=1,
xmax=66,
ytick={0, 0.25, 0.5, 0.75,1},
yticklabels={$0$, $25$, $50$, $75$, $100$},
xlabel style={font=\color{black}},
xlabel={\small Time (hrs)},
ymin=-0.05,
ymax=1.05,
yticklabel style = {font=\footnotesize,xshift=0ex},
xticklabel style = {font=\footnotesize,yshift=0ex},
ylabel style={font=\color{black}},
axis background/.style={fill=white},
xmajorgrids,
ymajorgrids,
legend style={legend cell align=left, align=left, draw=black}
]
\addplot [color=mycolor1,line width=1.5pt]
  table[row sep=crcr]{%
1	0.0311568440472172\\
2	0.0790264432082068\\
3	0.000151183895894746\\
4	0.000624022107757582\\
5	0.000826754823961551\\
6	0.659348007125459\\
7	0.999512583970855\\
8	0.999955869355574\\
11	0.999786374331563\\
12	0.699833117205102\\
13	0.100589405978582\\
14	0.0780717484030191\\
15	0.0581496042150604\\
16	0.0465789902712004\\
17	0.0380670351385959\\
18	0.0316497994644038\\
20	0.0202625439069806\\
21	0.0118809321922839\\
22	0.00523883862040009\\
23	3.7232971976664e-05\\
25	0.000182692388747796\\
26	0.0220562986533253\\
27	0.00226834427242295\\
28	0.724914858694035\\
29	0.59080500445026\\
30	0.966188525560312\\
31	0.999500935343519\\
32	0.999955634582875\\
33	0.999813518497533\\
35	0.999946625322565\\
36	0.695541884489046\\
37	0.131222007822714\\
38	0.10359535351779\\
39	0.0772085071459401\\
40	0.0618231024973568\\
41	0.0521727121000026\\
42	0.0423174148454706\\
43	0.0344256324024599\\
44	0.0289422454689401\\
45	0.0205766838294181\\
46	0.0146851642040389\\
47	0.010645152266477\\
48	2.43311056777884e-05\\
49	5.32037662992479e-05\\
50	0.179080400105448\\
51	0.241275379606734\\
52	0.125593853425499\\
53	0.00167185362775513\\
54	0.802518857236933\\
55	0.973197131023696\\
56	0.999958864872767\\
58	0.999920635740153\\
59	0.897645381100787\\
60	0.621197631486709\\
61	0.0720320802631846\\
62	0.0543346614633151\\
63	0.040575138994015\\
64	0.0299422434390522\\
65	0.0185347104986846\\
66	0.015757851210779\\
67	0.00989489163862345\\
69	0.00629176746068083\\
70	0.00261210834040071\\
71	2.59932830175558e-05\\
72	3.09962008770981e-05\\
};

\addplot [color=mycolor2,line width=1.5pt]
  table[row sep=crcr]{%
1	0.00623667448557796\\
2	0.0158146866101418\\
3	3.91494924798508e-05\\
4	0.0456837675232293\\
5	0.202816336701304\\
6	0.56742167667079\\
7	0.88156582421702\\
8	0.999990195761484\\
9	0.999979424674891\\
10	0.885660882331408\\
11	0.714181452961313\\
12	0.460507802632549\\
13	0.0900504851751123\\
14	0.0732451516651054\\
15	0.0580739653743194\\
16	0.0480286998923702\\
17	0.0398362708851892\\
18	0.0330688552716651\\
20	0.0208009785168457\\
21	0.0128638919213984\\
22	0.00606644967555781\\
23	0.000449837968133693\\
24	7.08516711256379e-06\\
25	5.90213419968677e-05\\
26	0.0461900622810134\\
27	0.13777795556517\\
28	0.393090171340091\\
29	0.524352853066915\\
30	0.780501802184261\\
31	0.955192137277265\\
32	0.999990910006261\\
34	0.99993751273675\\
35	0.858404865975146\\
36	0.53619341781949\\
37	0.130601294699062\\
38	0.106969103905342\\
39	0.0842600852083564\\
40	0.0694901649109738\\
41	0.0589642545251792\\
42	0.0486746911913798\\
43	0.0398866403722877\\
44	0.0328323761303153\\
45	0.0234676717060012\\
46	0.0159774741573386\\
47	0.0105211116006529\\
48	2.03517148378296e-05\\
49	1.63041310088374e-05\\
50	0.061016659234582\\
51	0.140601317166414\\
52	0.266302354091877\\
53	0.422796343347983\\
54	0.682501841727444\\
55	0.871630455136383\\
56	0.999988247662074\\
57	0.997995153548644\\
58	0.859280032262987\\
59	0.677833767514329\\
60	0.47679700388035\\
61	0.0840688137623857\\
62	0.0677993585650825\\
63	0.05408306065614\\
64	0.0431032468507908\\
65	0.0315580358404191\\
66	0.0265255354379264\\
67	0.0193622252213004\\
69	0.0106517608928414\\
70	0.00486187478020383\\
71	1.74910988022248e-05\\
72	6.20133832285319e-06\\
};

\addplot [color=mycolor3,line width=1.5pt]
  table[row sep=crcr]{%
1	0.010040044079858\\
2	0.0268908817179181\\
3	0.0531073454858273\\
4	0.125424090383873\\
5	0.258588500352161\\
6	0.455444620571754\\
7	0.640102467605473\\
8	0.771274650853343\\
9	0.754198223471022\\
10	0.611244117420853\\
11	0.471128679147654\\
12	0.300334715313539\\
13	0.0644719201639816\\
14	0.0533741133734509\\
15	0.0433573527424471\\
16	0.0365227699089985\\
17	0.0308277117186151\\
18	0.0260365771744091\\
20	0.0172537098945469\\
21	0.0116662810390267\\
22	0.00683098871664356\\
23	0.00279913975361978\\
24	3.13488949075236e-05\\
25	0.0137099716986171\\
26	0.0613206486408586\\
27	0.163823021330501\\
28	0.314754312230775\\
29	0.452507301986046\\
30	0.646052802653685\\
31	0.806100676159446\\
32	0.978906777327396\\
33	0.940807928722862\\
34	0.863497345775144\\
35	0.708719831513108\\
36	0.47710847689028\\
37	0.201877631973261\\
38	0.18538421164871\\
39	0.169496506438804\\
40	0.158970947004462\\
42	0.143829828213512\\
43	0.137378726474765\\
44	0.132051982903022\\
45	0.12513758697537\\
46	0.119510698190496\\
47	0.115278102428192\\
48	0.107667562357719\\
49	0.106376297809973\\
50	0.131357952686045\\
51	0.170221506940607\\
52	0.238978025545777\\
53	0.323700991558525\\
54	0.425948634163646\\
55	0.513453009784854\\
56	0.575833341254864\\
57	0.545152535013784\\
58	0.452406505250664\\
59	0.351698593198236\\
60	0.248327430890043\\
61	0.0461159899805637\\
62	0.0374397567754414\\
63	0.0300863722723506\\
64	0.0241806436284122\\
65	0.0179722322589413\\
66	0.0150954288685909\\
67	0.0111269239470744\\
69	0.00606859653063907\\
70	0.00278980503030368\\
71	1.41517262335356e-05\\
72	2.06842729255641e-06\\
};

\addplot [color=mycolor4,line width=1.5pt]
  table[row sep=crcr]{%
1	0.000502019358762595\\
2	0.00134458888358324\\
3	0.00265548173635466\\
4	0.00627137050123849\\
5	0.0129296402186441\\
6	0.0227725027704366\\
7	0.0320054432839783\\
8	0.0385640952671764\\
9	0.0377103691741922\\
11	0.0235569496934716\\
12	0.0150171561847685\\
13	0.00322375441922418\\
15	0.0021679926206275\\
18	0.00130191966867699\\
24	1.57281245094509e-06\\
25	0.000685606556615426\\
26	0.00306627505092649\\
27	0.00819171500279481\\
28	0.0157385471042062\\
29	0.0226265856946384\\
30	0.0323043000995398\\
31	0.040307085541528\\
32	0.048947956550947\\
33	0.047044297374768\\
34	0.0431794275293527\\
35	0.0354409509582752\\
36	0.0238605087090775\\
37	0.0100989433811662\\
39	0.00847989222513945\\
41	0.00757039597409914\\
45	0.00626190909568436\\
49	0.00532376537704238\\
50	0.0065728358351862\\
51	0.00851598164484813\\
52	0.0119537401782281\\
53	0.0161898266799625\\
54	0.0213021665726529\\
55	0.025677314763243\\
56	0.0287962620957956\\
57	0.0272617587344968\\
58	0.0226237149757367\\
60	0.01241829415045\\
61	0.00230620886264887\\
63	0.00150460537881258\\
66	0.00075492860825932\\
70	0.00013952090660041\\
72	1.03421669450654e-07\\
};

\end{axis}
\end{tikzpicture}%

%% file: Cobweb.tex
\tikzset{middlearrow/.style={
        decoration={markings,
            mark= at position 0.8 with {\arrow{#1}} ,
        },
        postaction={decorate}
    }
}
\begin{tikzpicture}

\definecolor{mycolor1}{HTML}{BE9063}%
\definecolor{mycolor10}{RGB}{136,176,197}
\definecolor{mycolor2}{RGB}{168,50,45}

\draw[->,] (-0.5,0) node[](a){}--(6,0) node[right]{$q$};
\draw[->,] (0,-0.5)--(0,5) node[above]{$\pi$};

\draw[mycolor1, line width=1pt, -] (0.4, 2)--(5.5,4) node[above]{\scriptsize {\color{black} $g_v(q)$}};
\draw[mycolor10, line width=1pt, -] (1.27, 4.5) node[above]{\scriptsize {\color{black} $g_k(q)$}}--(4,0.4) ;

\draw[mycolor2, line width=1pt, middlearrow={>}] (0.4, 0.6)--(3.84,0.6) ;
\draw[mycolor2, line width=1pt, middlearrow={>}] (3.87, 0.6)--(3.87,3.32) ;
\draw[mycolor2, line width=1pt, middlearrow={>}] (3.87, 3.35)--(2.06, 3.35) ;
\draw[mycolor2, line width=1pt, middlearrow={>}] (2.03, 3.35)--(2.03, 2.68) ;

\draw [mycolor10,fill=mycolor10] (3.87,0.6) circle [radius=0.075];
\draw [mycolor1,fill=mycolor1] (3.87,3.35) circle [radius=0.075];
\draw [mycolor10,fill=mycolor10] (2.03,3.35) circle [radius=0.075];
\draw [mycolor1,fill=mycolor1] (2.03,2.65) circle [radius=0.075];

\node at ($(a)+ (5.2,0.7)$) {\scriptsize ($q^{(1)}$, $\pi^{(0)}$)};
\node at ($(a)+ (5.2,3.3)$) {\scriptsize ($q^{(1)}$, $\pi^{(1)}$)};
\node at ($(a)+ (1.7,3.3)$){\scriptsize ($q^{(2)}$, $\pi^{(1)}$)};
\node at ($(a)+ (2.55,2.3)$)  {\scriptsize ($q^{(2)}$, $\pi^{(2)}$)};
\node at ($(a)+ (4,4.5)$)  {\small $\big|\!\frac{d g_k}{d q}\!\big|${{\tiny \textgreater}}$\!\frac{d g_v}{d q}$};

\end{tikzpicture}

%% file: probnoncon.tex
%
\definecolor{mycolor1}{HTML}{BE9063}%
\definecolor{mycolor2}{RGB}{168,50,45}
\definecolor{mycolor3}{HTML}{525B56}
\definecolor{mycolor4}{HTML}{A4978E}%
\definecolor{mycolor5}{rgb}{0.59400,0.18400,0.35600}%
\definecolor{mycolor6}{HTML}{253F5B}%
\definecolor{mycolor7}{HTML}{818A6F}%
\definecolor{mycolor8}{HTML}{D59B2D}
\definecolor{mycolor9}{RGB}{31,140,24}
\definecolor{mycolor10}{RGB}{136,176,197}

\tikzset{middlearrow/.style={
        decoration={markings,
            mark= at position 0.8 with {\arrow{#1}} ,
        },
        postaction={decorate}
    }
}

\begin{tikzpicture}[spy using outlines={circle, magnification=5.5,connect spies}]

\begin{axis}[%
width=2.9in,
height=1.6in,
at={(1.739in,0.849in)},
scale only axis,
xmin=0,
xmax=30,
domain=-5:5,enlargelimits=false,
xlabel style={font=\color{black}},
xlabel={\small Quantity traded (kWh)},
ymin=0.055,
ymax=.26,
ylabel style={font=\color{black}},
ylabel={\small Price $\pi$ (\$)},
ymajorgrids,
xmajorgrids,
yticklabel style = {font=\footnotesize,xshift=0ex},
xticklabel style = {font=\footnotesize,yshift=0ex},
axis background/.style={fill=white},
legend style={legend cell align=left, align=left, draw=black}
]
\addplot [color=mycolor2, line width=1.25pt]
  table[row sep=crcr]{%
0	0.0621408579846623\\
4.99999999350289	0.0621408579846623\\
4.99999999350289	0.0691261255616915\\
9.99999998512299	0.0691261255616915\\
9.99999998512299	0.0778983889384186\\
14.9999999379847	0.0778983889384186\\
14.9999999379847	0.0892579728246403\\
19.99999991982	0.0892579728246403\\
19.99999991982	0.104478009950299\\
24.9999997520808	0.104478009950299\\
};

\addplot [color=mycolor2, line width=0.25pt]
  table[row sep=crcr]{%
24.9999997520808	0.104478009950299\\
24.9999997520808	0.126041714904666\\
27.7599252488944	0.126041714904666\\
27.7599252488944	0.142256683424147\\
22.7599254098791	0.142256683424147\\
22.7599254098791	0.115373296335434\\
25.2599251389306	0.115373296335434\\
25.2599251389306	0.127408738688693\\
26.5099245344893	0.127408738688693\\
26.5099245344893	0.134421906192159\\
25.2599246427683	0.134421906192159\\
25.2599246427683	0.12740873605058\\
25.8849209549673	0.12740873605058\\
25.8849209549673	0.130820930223411\\
25.5724211763055	0.130820930223411\\
25.5724211763055	0.129092122811286\\
25.7286697268963	0.129092122811286\\
25.7286697268963	0.129950744519675\\
25.6505457810398	0.129950744519675\\
25.6505457810398	0.129520004663661\\
25.6896034370628	0.129520004663661\\
25.6896034370628	0.129734991838482\\
25.6700846976492	0.129734991838482\\
25.6700846976492	0.129627464289349\\
25.6798467142586	0.129627464289349\\
25.6798467142586	0.129681220267972\\
};

\addplot [color=mycolor1, line width=2pt]
  table[row sep=crcr]{%
0	0.244941262456212\\
0.280403285342366	0.241578785493232\\
0.560806570684736	0.238351031814325\\
0.841209856027103	0.235249239431962\\
1.12161314136947	0.232265423917987\\
1.40201642671184	0.229392291879432\\
1.68241971205421	0.226623165918166\\
1.96282299739658	0.223951919315581\\
2.24322628273895	0.221372918986994\\
2.52362956808131	0.218880975496248\\
2.80403285342368	0.216471299120524\\
3.08443613876605	0.214139461118545\\
3.36483942410841	0.211881359489158\\
3.64524270945078	0.209693188617848\\
3.92564599479315	0.207571412299902\\
4.20604928013552	0.205512739705174\\
4.48645256547789	0.203514103912696\\
5.04725913616262	0.199685681289559\\
5.60806570684736	0.196065404692082\\
6.1688722775321	0.192635073372198\\
6.72967884821683	0.189378639431872\\
7.29048541890157	0.1862818911146\\
7.8512919895863	0.183332191241011\\
8.41209856027104	0.180518259789075\\
8.97290513095577	0.177829992075996\\
9.53371170164051	0.175258305851074\\
10.0945182723252	0.172795012019193\\
10.65532484301	0.170432704796809\\
11.2161314136947	0.168164667940051\\
11.7769379843795	0.165984794337778\\
12.3377445550642	0.163887516775109\\
13.1789544110913	0.160885458150837\\
14.0201642671184	0.158042483869799\\
14.8613741231455	0.155345023543049\\
15.7025839791726	0.152781073961169\\
16.5437938351997	0.150339973825172\\
17.3850036912268	0.14801221682966\\
18.2262135472539	0.145789295613135\\
19.067423403281	0.143663570726659\\
19.9086332593081	0.141628160014772\\
21.0302464006776	0.139044104771596\\
22.1518595420471	0.136596460546485\\
};

\addplot [color=mycolor1, line width=0.25pt]
  table[row sep=crcr]{%
22.1518595420471	0.136596460546485\\
23.2734726834165	0.134273633959065\\
24.395085824786	0.132065365718617\\
25.5166989661555	0.129962539521831\\
26.638312107525	0.127957023377071\\
27.7599252488944	0.126041537046969\\
};

\addplot [color=mycolor10, line width=2pt]
  table[row sep=crcr]{%
0	0.0621466629365166\\
1.40201642671184	0.0639602875487846\\
2.80403285342368	0.0658834679393152\\
4.20604928013552	0.067926456898757\\
5.60806570684736	0.0701008288682274\\
6.72967884821683	0.0719436780166127\\
7.8512919895863	0.073886506399969\\
8.97290513095577	0.0759376905928661\\
10.0945182723252	0.0781065706085045\\
11.2161314136947	0.0804035925860269\\
12.3377445550642	0.0828404776198219\\
13.1789544110913	0.0847678032342962\\
14.0201642671184	0.0867873779760338\\
14.8613741231455	0.0889059981964628\\
15.7025839791726	0.0911311452347263\\
16.5437938351997	0.0934710739907842\\
17.3850036912268	0.09593491561197\\
18.2262135472539	0.0985327969888345\\
19.067423403281	0.101275980361418\\
19.9086332593081	0.10417702709946\\
20.4694398299929	0.106205632092145\\
21.0302464006776	0.108315187483694\\
21.5910529713623	0.110510645171225\\
22.1518595420471	0.112797369681534\\
};

\addplot [color=mycolor10, line width=0.25pt]
  table[row sep=crcr]{%
22.1518595420471	0.112797369681534\\
22.7126661127318	0.115181182086491\\
23.2734726834165	0.117668409650928\\
23.8342792541013	0.120265942105679\\
24.395085824786	0.122981295600937\\
24.9558923954707	0.125822685591583\\
25.5166989661555	0.128799110144996\\
26.0775055368402	0.131920445452995\\
26.638312107525	0.13519755568646\\
27.1991186782097	0.138642419770573\\
27.7599252488944	0.142268278202131\\
};

\addplot[only marks, mark=*, mark options={solid, color=mycolor1}, mark size=2.5pt, draw=mycolor1] table[row sep=crcr]{%
x	y\\
4.99999999350289	0.0621408579846614\\
9.99999998512299	0.0691261255616928\\
14.9999999379847	0.0778983889384176\\
19.99999991982	0.0892579728246413\\
};

\addplot[only marks, mark=*, mark options={solid, color=mycolor10, color=mycolor1}, mark size=0.4pt, draw=mycolor1] table[row sep=crcr]{%
x	y\\
24.9999997520808	0.104478009950299\\
27.7599252488944	0.126041714904666\\
22.7599254098791	0.142256683424148\\
25.2599251389306	0.115373296335436\\
26.5099245344893	0.127408738688694\\
25.2599246427683	0.134421906192159\\
25.8849209549673	0.127408736050578\\
25.5724211763055	0.130820930223412\\
25.7286697268963	0.129092122811285\\
25.6505457810398	0.129950744519675\\
25.6896034370628	0.129520004663662\\
25.6700846976492	0.129734991838483\\
25.6798467142586	0.129627464289348\\
};

\addplot[only marks, mark=*, mark options={solid, color=mycolor10}, mark size=2.5pt, draw=mycolor10] table[row sep=crcr]{%
x	y\\
0	0.0621408579846614\\
4.99999999350289	0.0691261255616928\\
9.99999998512299	0.0778983889384176\\
14.9999999379847	0.0892579728246413\\
19.99999991982	0.104478009950299\\
};

\addplot[only marks, mark=*, mark options={solid, color=mycolor10}, mark size=0.4pt, draw=mycolor10] table[row sep=crcr]{%
x	y\\
24.9999997520808	0.126041714904666\\
27.7599252488944	0.142256683424148\\
22.7599254098791	0.115373296335436\\
25.2599251389306	0.127408738688694\\
26.5099245344893	0.134421906192159\\
25.2599246427683	0.127408736050578\\
25.8849209549673	0.130820930223412\\
25.5724211763055	0.129092122811285\\
25.7286697268963	0.129950744519675\\
25.6505457810398	0.129520004663662\\
25.6896034370628	0.129734991838483\\
25.6700846976492	0.129627464289348\\
25.6798467142586	0.129681220267972\\
};

\addplot[only marks, mark=asterisk, mark options={}, mark size=0.5pt, draw=mycolor2] table[row sep=crcr]{%
x	y\\
25.6764073196101	0.129668566207077\\
};

\coordinate (spypoint) at (axis cs:25.85,0.128);
\coordinate (spyviewer) at (axis cs:21,0.215);
\spy[width=1in,height=0.6in] on (spypoint) in node [fill=white] at (spyviewer);

\draw[mycolor2,line width=0.25pt,middlearrow={>}] (axis cs: 24.9999997520808,	0.104478009950299)--(axis cs: 24.9999997520808,	0.126041714904666);
\draw[mycolor2, line width=0.25pt,middlearrow={>}] (axis cs: 27.7599252488944,	0.126041714904666)--(axis cs: 27.7599252488944,	0.142256683424147);
\draw[mycolor2, line width=0.25pt,middlearrow={>}] (axis cs: 27.7599252488944,	0.142256683424147)--(axis cs: 22.7599254098791,	0.142256683424147);
\draw[mycolor2, line width=0.25pt,middlearrow={>}] (axis cs: 22.7599254098791,	0.142256683424147)--(axis cs: 22.7599254098791,	0.115373296335434);
\draw[mycolor2, line width=0.25pt,middlearrow={>}] (axis cs: 22.7599254098791,	0.115373296335434)--(axis cs: 25.2599251389306,	0.115373296335434);
\draw[mycolor2, line width=0.25pt,middlearrow={>}] (axis cs: 25.2599251389306,	0.115373296335434)--(axis cs: 25.2599251389306,	0.127408738688693);
\draw[mycolor2, line width=0.01pt,middlearrow={>}] (axis cs: 26.5099245344893,	0.127408738688693)--(axis cs: 26.5099245344893,	0.134421906192159);
\draw[mycolor2, line width=0.1pt,middlearrow={>}] (axis cs: 26.5099245344893,	0.134421906192159)--(axis cs: 25.2599246427683,	0.134421906192159);
\draw[mycolor2, line width=0.1pt,middlearrow={>}] (axis cs: 25.2599246427683,	0.134421906192159)--(axis cs: 25.2599246427683,	0.12740873605058);

\draw[mycolor2, line width=1.5pt,middlearrow={>}] (axis cs: 0,	0.0621408579846623)--(axis cs: 4.99999999350289,	0.0621408579846623);
\draw[mycolor2, line width=1.5pt,middlearrow={>}] (axis cs: 4.99999999350289,	0.0691261255616915)--(axis cs: 9.99999998512299,	0.0691261255616915);
\draw[mycolor2, line width=1.5pt,middlearrow={>}] (axis cs:  9.99999998512299,	0.0778983889384186)--(axis cs: 14.9999999379847,	0.0778983889384186);
\draw[mycolor2, line width=1.5pt,middlearrow={>}] (axis cs:  14.9999999379847,	0.0892579728246403)--(axis cs: 19.99999991982,	0.0892579728246403);
\draw[mycolor2, line width=1.5pt,middlearrow={>}] (axis cs:  19.99999991982,	0.104478009950299)--(axis cs: 24.9999997520808,	0.104478009950299);

\end{axis}
\end{tikzpicture}%

%% file: parameter_experiment.tex
%
\definecolor{mycolor3}{HTML}{D59B2D}
\definecolor{mycolor2}{RGB}{31,140,24}
\definecolor{mycolor1}{RGB}{136,176,197}
\definecolor{mycolor2}{RGB}{168,50,45}

\begin{tikzpicture}

\begin{axis}[%
width=2.9in,
height=1.25in,
at={(1.739in,2.629in)},
scale only axis,
xmin=0,
xmax=1,
xlabel style={font=\color{black}},
xlabel={\small $\gamma$},
ylabel={\small $\#$ Iterations},
ymode=log,
ymin=10,
ymax=1000,
xmajorgrids,
yticklabel style = {font=\footnotesize,xshift=0ex},
xticklabel style = {font=\footnotesize,yshift=0ex},
ymajorgrids,
yminorticks=true,
axis background/.style={fill=white},
legend style={legend cell align=left, align=left, draw=black, font=\small, draw=none, legend columns=-1, at={(1.04,1.25)}}
]
\addplot [color=mycolor1, line width=2.0pt]
  table[row sep=crcr]{%
0.05	98.265\\
0.1	41.1516666666667\\
0.15	23.5216666666667\\
0.2	18.19\\
0.25	15.3816666666667\\
0.3	13.9533333333333\\
0.35	13.3483333333333\\
0.4	13.2033333333333\\
0.45	13.28\\
0.5	13.69\\
0.55	14.21\\
0.6	14.92\\
0.65	16.1016666666667\\
0.7	17.5433333333333\\
0.75	19.7466666666667\\
0.8	22.94\\
0.85	28.025\\
0.9	38.7266666666667\\
0.95	69.42\\
};
\addlegendentry{$\delta^{(0)} \in [0.1, 0.7$)~}

\addplot [color=mycolor2, line width=2.0pt]
  table[row sep=crcr]{%
0.05	102.212857142857\\
0.1	36.7571428571429\\
0.15	19.7385714285714\\
0.2	15.1271428571429\\
0.25	12.9771428571429\\
0.3	12.2128571428571\\
0.35	11.9414285714286\\
0.4	12.0742857142857\\
0.45	12.5414285714286\\
0.5	13.1071428571429\\
0.55	14.0285714285714\\
0.6	15.11\\
0.65	16.6328571428571\\
0.7	18.6271428571429\\
0.75	21.1942857142857\\
0.8	24.8871428571429\\
0.85	31.1057142857143\\
0.9	43.4171428571429\\
0.95	79.6485714285714\\
};
\addlegendentry{$\delta^{(0)} \in [0.7, 1.3$)~}

\addplot [color=mycolor3, line width=2.0pt]
  table[row sep=crcr]{%
0.05	94.1714285714286\\
0.1	35.1557142857143\\
0.15	19.8128571428571\\
0.2	15.6642857142857\\
0.25	13.3\\
0.3	12.4814285714286\\
0.35	12.2914285714286\\
0.4	12.5828571428571\\
0.45	13.0371428571429\\
0.5	13.75\\
0.55	14.5242857142857\\
0.6	15.7471428571429\\
0.65	17.2671428571429\\
0.7	19.43\\
0.75	22.2385714285714\\
0.8	26.1542857142857\\
0.85	32.7671428571429\\
0.9	45.8685714285714\\
0.95	84.9228571428571\\
};
\addlegendentry{$\delta^{(0)} \in [1.3, 2$]~}

\addplot [color=mycolor1, dashed, line width=1.75pt, forget plot]
  table[row sep=crcr]{%
0.05	711\\
0.1	233\\
0.15	157\\
0.2	119\\
0.25	96\\
0.3	80\\
0.35	70\\
0.4	61\\
0.45	56\\
0.5	52\\
0.55	49\\
0.6	47\\
0.65	45\\
0.7	45\\
0.75	46\\
0.8	55\\
0.85	75\\
0.9	114\\
0.95	231\\
};
\addplot [color=mycolor2, dashed, line width=1.75pt, forget plot]
  table[row sep=crcr]{%
0.05	706\\
0.1	168\\
0.15	101\\
0.2	64\\
0.25	39\\
0.3	34\\
0.35	33\\
0.4	26\\
0.45	29\\
0.5	28\\
0.55	29\\
0.6	31\\
0.65	35\\
0.7	41\\
0.75	50\\
0.8	59\\
0.85	79\\
0.9	121\\
0.95	246\\
};
\addplot [color=mycolor3, dashed, line width=1.75pt, forget plot]
  table[row sep=crcr]{%
0.05	587\\
0.1	166\\
0.15	90\\
0.2	55\\
0.25	45\\
0.3	37\\
0.35	28\\
0.4	29\\
0.45	28\\
0.5	29\\
0.55	30\\
0.6	34\\
0.65	39\\
0.7	43\\
0.75	53\\
0.8	64\\
0.85	83\\
0.9	125\\
0.95	255\\
};
\end{axis}

\end{tikzpicture}%

%% file: test_boxplot.tex
%
%
\definecolor{mycolor1}{HTML}{BE9063}%
\definecolor{mycolor2}{RGB}{168,50,45}
\definecolor{mycolor3}{HTML}{525B56}
\definecolor{mycolor4}{HTML}{A4978E}%
\definecolor{mycolor5}{rgb}{0.59400,0.18400,0.35600}%
\definecolor{mycolor6}{HTML}{253F5B}%
\definecolor{mycolor7}{HTML}{818A6F}
\definecolor{mycolor8}{HTML}{D59B2D}
\definecolor{mycolor9}{RGB}{31,140,24}
\definecolor{mycolor10}{RGB}{136,176,197}

\begin{tikzpicture}

\begin{axis}[%
width=2.9in,
height=1.25in,
at={(1.739in,0.849in)},
scale only axis,
xlabel style={font=\color{black}},
xlabel={\small Total Battery Capacity [kWh]},
ylabel style={font=\color{black}},
ylabel={\small \# Iterations},
xmin=0.5,
xmax=5.5,
xtick={1,2,3,4,5},
xticklabels={$15$,$25$,$40$,$80$,$300$},
ymin=0,
ymax=450,
ytick={0, 100, 200, 300, 400},
yticklabels={$0$,$100$,$200$, $300$, $400$},
ymajorgrids,
xmajorgrids,
yticklabel style = {font=\footnotesize,xshift=0ex},
xticklabel style = {font=\footnotesize,yshift=0ex},
axis background/.style={fill=white}
]
\addplot [color=mycolor2, dashed, forget plot, line width=1.05pt]
  table[row sep=crcr]{%
1	88\\
1	178\\
};
\addplot [color=mycolor2, dashed, forget plot, line width=1.05pt]
  table[row sep=crcr]{%
2	88\\
2	168\\
};
\addplot [color=mycolor2, dashed, forget plot, line width=1.05pt]
  table[row sep=crcr]{%
3	106.5\\
3	204\\
};
\addplot [color=mycolor2, dashed, forget plot, line width=1.05pt]
  table[row sep=crcr]{%
4	138.5\\
4	244\\
};
\addplot [color=mycolor2, dashed, forget plot, line width=1.05pt]
  table[row sep=crcr]{%
5	197\\
5	431\\
};
\addplot [color=mycolor2, dashed, forget plot, line width=1.05pt]
  table[row sep=crcr]{%
1	16\\
1	25\\
};
\addplot [color=mycolor2, dashed, forget plot, line width=1.05pt]
  table[row sep=crcr]{%
2	14\\
2	31\\
};
\addplot [color=mycolor2, dashed, forget plot, line width=1.05pt]
  table[row sep=crcr]{%
3	16\\
3	37\\
};
\addplot [color=mycolor2, dashed, forget plot, line width=1.05pt]
  table[row sep=crcr]{%
4	14\\
4	38\\
};
\addplot [color=mycolor2, dashed, forget plot, line width=1.05pt]
  table[row sep=crcr]{%
5	21\\
5	40\\
};
\addplot [color=mycolor2, forget plot, line width=1.05pt]
  table[row sep=crcr]{%
0.875	178\\
1.125	178\\
};
\addplot [color=mycolor2, forget plot, line width=1.05pt]
  table[row sep=crcr]{%
1.875	168\\
2.125	168\\
};
\addplot [color=mycolor2, forget plot, line width=1.05pt]
  table[row sep=crcr]{%
2.875	204\\
3.125	204\\
};
\addplot [color=mycolor2, forget plot, line width=1.05pt]
  table[row sep=crcr]{%
3.875	244\\
4.125	244\\
};
\addplot [color=mycolor2, forget plot, line width=1.05pt]
  table[row sep=crcr]{%
4.875	431\\
5.125	431\\
};
\addplot [color=mycolor2, forget plot, line width=1.05pt]
  table[row sep=crcr]{%
0.875	16\\
1.125	16\\
};
\addplot [color=mycolor2, forget plot, line width=1.05pt]
  table[row sep=crcr]{%
1.875	14\\
2.125	14\\
};
\addplot [color=mycolor2, forget plot, line width=1.05pt]
  table[row sep=crcr]{%
2.875	16\\
3.125	16\\
};
\addplot [color=mycolor2, forget plot, line width=1.05pt]
  table[row sep=crcr]{%
3.875	14\\
4.125	14\\
};
\addplot [color=mycolor2, forget plot, line width=1.05pt]
  table[row sep=crcr]{%
4.875	21\\
5.125	21\\
};
\addplot [color=mycolor10, forget plot, line width=1.05pt]
  table[row sep=crcr]{%
0.75	25\\
0.75	88\\
1.25	88\\
1.25	25\\
0.75	25\\
};
\addplot [color=mycolor10, forget plot, line width=1.05pt]
  table[row sep=crcr]{%
1.75	31\\
1.75	88\\
2.25	88\\
2.25	31\\
1.75	31\\
};
\addplot [color=mycolor10, forget plot, line width=1.05pt]
  table[row sep=crcr]{%
2.75	37\\
2.75	106.5\\
3.25	106.5\\
3.25	37\\
2.75	37\\
};
\addplot [color=mycolor10, forget plot, line width=1.05pt]
  table[row sep=crcr]{%
3.75	38\\
3.75	138.5\\
4.25	138.5\\
4.25	38\\
3.75	38\\
};
\addplot [color=mycolor10, forget plot, line width=1.05pt]
  table[row sep=crcr]{%
4.75	40\\
4.75	197\\
5.25	197\\
5.25	40\\
4.75	40\\
};
\addplot [color=mycolor1, forget plot, line width=1.05pt]
  table[row sep=crcr]{%
0.75	54\\
1.25	54\\
};
\addplot [color=mycolor1, forget plot, line width=1.05pt]
  table[row sep=crcr]{%
1.75	54\\
2.25	54\\
};
\addplot [color=mycolor1, forget plot, line width=1.05pt]
  table[row sep=crcr]{%
2.75	61\\
3.25	61\\
};
\addplot [color=mycolor1, forget plot, line width=1.05pt]
  table[row sep=crcr]{%
3.75	63\\
4.25	63\\
};
\addplot [color=mycolor1, forget plot, line width=1.05pt]
  table[row sep=crcr]{%
4.75	88\\
5.25	88\\
};
\end{axis}
\end{tikzpicture}%

%% file: bidding_plots.tex
%
\definecolor{mycolor2}{HTML}{F1231B}%
\definecolor{mycolor3}{HTML}{D59B2D}
\definecolor{mycolor1}{RGB}{136,176,197}
\definecolor{mycolor4}{HTML}{253F5B}%
\definecolor{mycolor5}{HTML}{818A6F}
\definecolor{mycolor6}{RGB}{168,50,45}
\begin{tikzpicture}

\begin{axis}[%
width=2.9in,
height=1.25in,
at={(1.739in,2.849in)},
scale only axis,
xmin=1,
xmax=24,
ymin=0,
ymax=18,
xtick={2,4,6,8,10,12,14,16,18,20,22,24},
xticklabels={$2$,$4$,$6$,$8$,$10$,$12$,$14$,$16$,$18$,$20$,$22$,$24$},
ylabel style={font=\color{black}},
ylabel={\small Price [\$/kWh]},
yticklabel style = {font=\footnotesize,xshift=0ex},
xticklabel style = {font=\footnotesize,yshift=0ex},
axis background/.style={fill=white},
xmajorgrids,
ymajorgrids,
legend style={legend cell align=left, align=left, draw=black, draw=none, legend columns=-1, at={(.88,1.2)}}
]
\addplot [color=mycolor1, line width=1.5pt, mark=*, mark size=1.5pt]
  table[row sep=crcr]{%
1	5.37389377700304\\
2	5.37389378945041\\
3	5.37389458332149\\
4	5.37389335139887\\
5	5.37389206626265\\
6	5.37389095027010\\
7	5.37389539671740\\
8	5.37389363668053\\
9	5.37389784380334\\
10	5.37389442381269\\
11	5.37389420785485\\
12	5.37389041223069\\
13	5.37388845063122\\
14	5.37389306430259\\
15	5.37389352653150\\
16	5.37389356523825\\
17	5.37389582070952\\
18	5.37389479773021\\
19	5.37389849141591\\
20	5.37389970205574\\
21	5.37389407894658\\
22	5.37389632685958\\
23	5.37389490017657\\
24	5.37389627548336\\
};
\addlegendentry{${\pi}_\text{centr}~$}

\addplot [color=mycolor2, dashed, line width=1.5pt, mark=*, mark size=1.5pt, mark options={solid}]
  table[row sep=crcr]{%
1	2.14637514129499\\
2	2.14636677205890\\
3	2.14636537973439\\
4	2.14636558267531\\
5	2.14637787311810\\
6	2.53250996395946\\
7	10.7106769215316\\
8	10.7106769215316\\
9	10.7106769215316\\
10	10.7106769215316\\
11	10.7106769215316\\
12	10.7106769215316\\
13	13.4075630363077\\
14	13.4075630363077\\
15	13.4075630363077\\
16	13.4075630363077\\
17	13.4075630363077\\
18	13.4075630363077\\
19	13.4075630363077\\
20	13.4075630363077\\
21	3.36330049186945\\
22	3.36330552324653\\
23	3.36330384612083\\
24	2.22595558874309\\
};
\addlegendentry{${\pi}_1~$}

\addplot [color=mycolor3, line width=1.5pt, mark=*, mark size=1.5pt, mark options={solid}]
  table[row sep=crcr]{%
1	2.25957726354216\\
2	2.25957725343184\\
3	2.25957844747441\\
4	2.25957787998281\\
5	2.25957586459774\\
6	2.54185762334186\\
7	9.41085383202483\\
8	9.41085388282686\\
9	9.41085392689820\\
10	9.41085392750266\\
11	9.41085396286885\\
12	9.41085395573085\\
13	11.9385413483740\\
14	11.9385402732401\\
15	11.9385393636183\\
16	11.9385397498497\\
17	11.9385398707103\\
18	11.9385399157455\\
19	11.9385399849559\\
20	11.9385400756179\\
21	11.9385401759760\\
22	11.9385402603085\\
23	11.9385403293127\\
24	2.55278593210375\\
};
\addlegendentry{${\pi}_2~$}

\addplot [color=mycolor4, dashed, line width=1.5pt, mark=*, mark size=1.5pt, mark options={solid}]
  table[row sep=crcr]{%
1	2.14722891607624\\
2	2.14722878758854\\
3	2.14722778135348\\
4	2.14722659862515\\
5	2.14722749247320\\
6	2.53250735408619\\
7	10.7109793802651\\
8	10.7109793944062\\
9	10.7109794053058\\
10	10.7109794054584\\
11	10.7109794132240\\
12	10.7109794112354\\
13	13.4153823885965\\
14	13.4153823901494\\
15	13.4153824115958\\
16	13.4153824357102\\
17	13.4153824456604\\
18	13.4153824570214\\
19	13.4153824747060\\
20	13.4153824979097\\
21	3.39694777682871\\
22	3.39694779946776\\
23	3.39694781634238\\
24	2.22879381881266\\
};
\addlegendentry{${\pi}_3~$}

\addplot [color=mycolor5, loosely dashed, line width=1.5pt, mark=*, mark size=1.5pt,mark options={solid}]
  table[row sep=crcr]{%
1	2.14795883049162\\
2	2.14795792320732\\
3	2.14795754187980\\
4	2.14795839411426\\
5	2.14795783976739\\
6	2.53252051174739\\
7	10.6396033433008\\
8	10.6396033420609\\
9	10.6396033411047\\
10	10.6396033410923\\
11	10.6396033404113\\
12	10.6396033405062\\
13	13.2617066850800\\
14	13.2617066848593\\
15	13.2617066829757\\
16	13.2617066808574\\
17	13.2617066799830\\
18	13.2617066789843\\
19	13.2617066774309\\
20	13.2617066753932\\
21	3.41770486579456\\
22	3.41770486393169\\
23	3.41770486220024\\
24	2.23145420817470\\
};
\addlegendentry{${\pi}_4~$}

\addplot [color=mycolor6, dotted, line width=1.5pt, mark=*, mark size=1.5pt, mark options={solid}]
  table[row sep=crcr]{%
1	2.14722891607624\\
2	2.14722878758854\\
3	2.14722778135348\\
4	2.14722659862515\\
5	2.14722749247320\\
6	2.53250735408619\\
7	10.7109793802651\\
8	10.7109793944062\\
9	10.7109794053058\\
10	10.7109794054584\\
11	10.7109794132240\\
12	10.7109794112354\\
13	13.4153823885965\\
14	13.4153823901494\\
15	13.4153824115958\\
16	13.4153824357102\\
17	13.4153824456604\\
18	13.4153824570214\\
19	13.4153824747060\\
20	13.4153824979097\\
21	3.39694777682871\\
22	3.39694779946776\\
23	3.39694781634238\\
24	2.22879381881266\\
};
\addlegendentry{${\pi}_5$}

\end{axis}

\end{tikzpicture}%